\documentclass[11pt,a4paper]{article}

\usepackage[T1]{fontenc}
\usepackage[utf8]{inputenc}

\usepackage[
    left=2.7cm,
    right=2.7cm,
    top=2.5cm,
    bottom=2.7cm
]{geometry}

\usepackage{mathtools,amssymb,amsthm}
\usepackage{stmaryrd,mathrsfs,calligra}
\usepackage{stix}
\usepackage{dsfont}
\usepackage{microtype}

\usepackage{graphicx,xcolor,float}
\graphicspath{{./images/}}

\usepackage[all]{xy}
\usepackage{tikz,tikz-cd}

\usepackage[
    backend=biber,
    style=alphabetic
]{biblatex}

\usepackage{titling}

\pretitle{\centering\LARGE\bfseries}
\posttitle{\par\vspace{0.8em}}

\preauthor{\centering\normalsize\scshape}
\postauthor{\par}

\predate{}
\postdate{}
\date{}

\usepackage[runin]{abstract}

\abslabeldelim{.\ }
\AtBeginDocument{\setlength{\absparindent}{0pt}}

\usepackage{titlesec}

\titleformat{\section}
  {\normalfont\large\bfseries}
  {\thesection.}{0.65em}{}

\titleformat{\subsection}
  {\normalfont\normalsize\bfseries}
  {\thesubsection.}{0.65em}{}

\titlespacing*{\section}
  {0pt}{2.0ex plus .5ex minus .2ex}{0.9ex plus .2ex}

\titlespacing*{\subsection}
  {0pt}{1.4ex plus .3ex minus .2ex}{0.6ex plus .2ex}

\usepackage[hidelinks]{hyperref}
\usepackage[nameinlink,noabbrev]{cleveref}

\theoremstyle{plain}
\newtheorem{theorem}{Theorem}[section]
\newtheorem{lemma}[theorem]{Lemma}
\newtheorem{proposition}[theorem]{Proposition}
\newtheorem{corollary}[theorem]{Corollary}

\theoremstyle{definition}
\newtheorem{definition}[theorem]{Definition}

\theoremstyle{remark}
\newtheorem{remark}[theorem]{Remark}

\titleformat{\section}
  {\normalfont\Large\bfseries}
  {\thesection.}{0.65em}{}

\titleformat{\subsection}
  {\normalfont\large\bfseries}
  {\thesubsection.}{0.65em}{}

\title{The Symplectic Nature of reducible noncentral Stratum in $\mathrm{SU}(2)$ Character Variety}
\author{Shuo Liu}
\date{}  

\begin{document}
\maketitle
\begin{abstract}
Let \(\Sigma_g\) be a closed oriented surface of genus \(g\geq2\). The purpose of this paper is to use an explicit example, the representation of the surface group $\pi_1(\Sigma_g)$ into $\mathrm{SU}(2)$, to investigate the deformation space of the representation variety at a type of singular points, reducible but noncentral representations, and explore the geometric structure of the corresponding stratum in the representation variety and the character variety. In particular, the reduced tangent quadratic cone encodes the symplectic structure at the reducible noncentral stratum of the character variety. Such a symplectic structure fits in the ambient character variety $\chi\big(\pi_1(\Sigma_g),\mathrm{SU}(2)\big)$, where the top stratum of smooth points admits the Narasimhan-Atiyah-Bott-Goldman symplectic structure.

\end{abstract}
\tableofcontents

\newpage

\section{Introduction}
Let $\Sigma_g$ be a closed, connected, oriented surface of genus
$g \geq 2$, and write
\[
\mathcal{R}
:= \mathrm{Hom}\bigl(\pi_1(\Sigma_g),\mathrm{SU}(2)\bigr).
\]
Consider the chain of subgroup inclusions
\[
\mathrm{O}(1)
\hookrightarrow \mathrm{U}(1)
\hookrightarrow \mathrm{SU}(2),
\]
where $\mathrm{O}(1)=\{\pm1\}$ and the second inclusion is given by
$z \mapsto \operatorname{diag}(z,z^{-1})$.
These inclusions induce embeddings
\[
\mathrm{Hom}\bigl(\pi_1(\Sigma_g),\mathrm{O}(1)\bigr)
\hookrightarrow
\mathrm{Hom}\bigl(\pi_1(\Sigma_g),\mathrm{U}(1)\bigr)
\hookrightarrow \mathcal{R}.
\]
Denote the images of the first and second representation spaces
in $\mathcal{R}$ by $\mathcal{Y}_1$ and $\mathcal{Y}$, respectively.
Thus, $\mathcal{Y}_1$ is the finite set of central representations,
and
\[
\mathcal{Y}_0 := \mathcal{Y}\setminus\mathcal{Y}_1
\]
is the open subset of noncentral representations in $\mathcal{Y}$.\\\\
Every reducible $\mathrm{SU}(2)$-representation is conjugate to a
representation with image in the fixed maximal torus
$\mathrm{U}(1)$. Consequently, the locus of noncentral reducible
representations is the conjugation saturation
\[
\widetilde{\mathcal{Y}}_0
:= \mathrm{Ad}_{\mathrm{SU}(2)}\mathcal{Y}_0.
\]
Finally, denote the irreducible locus by
\[
\mathcal{R}^{\mathrm{irr}}
:= \mathrm{Hom}\bigl(\pi_1(\Sigma_g),
                    \mathrm{SU}(2)\bigr)^{\mathrm{irr}}.
\]
These loci give the following stratification.

\begin{proposition}{(\ref{StratifiedRepresentationVariety})}
The representation space $\mathcal{R}$, viewed as a closed subset
of $\mathrm{SU}(2)^{2g}$ via the standard presentation of
$\pi_1(\Sigma_g)$, admits the orbit-type decomposition
\[
\mathcal{R}
=
\mathcal{Y}_1
\bigsqcup
\widetilde{\mathcal{Y}}_0
\bigsqcup
\mathcal{R}^{\mathrm{irr}}.
\]
The stabilizers under conjugation on these three pieces are
conjugate to $\mathrm{SU}(2)$, $\mathrm{U}(1)$, and
$\{\pm I\}$, respectively. Their closures in $\mathcal{R}$ satisfy
\[
\overline{\mathcal{Y}_1}=\mathcal{Y}_1,
\qquad
\overline{\widetilde{\mathcal{Y}}_0}
=
\widetilde{\mathcal{Y}}_0\bigsqcup\mathcal{Y}_1,
\qquad
\overline{\mathcal{R}^{\mathrm{irr}}}=\mathcal{R}.
\]
In particular, under the frontier order
$S\preceq T$ defined by $S\subseteq\overline{T}$, these pieces form
the chain
\[
\mathcal{Y}_1
\prec
\widetilde{\mathcal{Y}}_0
\prec
\mathcal{R}^{\mathrm{irr}}.
\]
\end{proposition}
The middle stratum $\widetilde{\mathcal{Y}}_0$ is the main object
of study in this paper. To investigate the local structure of
$\mathcal{R}$ along this stratum, we first recall the tangent
quadratic cone and the quadraticity theorem of
Goldman--Millson~\cite{BillMillsonDeforamtion}.

\begin{definition}[Tangent quadratic cone]
\label{def:tangent-quadratic-cone}
For $\rho_0\in\mathcal{R}$, the tangent quadratic cone to
$\mathcal{R}$ at $\rho_0$ is
\[
\mathcal{L}_{\rho_0}
:=
\left\{
u\in Z^1\bigl(\pi_1(\Sigma_g),
             \mathfrak{su}(2)_{\mathrm{Ad}_{\rho_0}}\bigr)
\;\middle|\;
[u,u]=0
\text{ in }
H^2\bigl(\pi_1(\Sigma_g),
        \mathfrak{su}(2)_{\mathrm{Ad}_{\rho_0}}\bigr)
\right\}.
\]
Here $[u,u]$ denotes the cohomology class of the cup product
of $u$ with itself, with the Lie bracket on $\mathfrak{su}(2)$
as the coefficient pairing. Thus, $\mathcal{L}_{\rho_0}$
consists of infinitesimal deformations whose second-order
obstruction vanishes.
\end{definition}
Recall that a real analytic space is \emph{quadratic at a point}
if its analytic germ there is isomorphic to the germ at the
vertex of a cone defined by homogeneous quadratic equations.

\begin{theorem}[Goldman--Millson, 1988]
\label{thm:goldman-millson-quadraticity}
The representation space
\[
\mathcal{R}
=
\mathrm{Hom}\bigl(\pi_1(\Sigma_g),\mathrm{SU}(2)\bigr)
\]
is quadratic at every representation $\rho_0\in\mathcal{R}$.
More precisely, there is an isomorphism of pointed real
analytic germs
\[
\big(\mathcal{R},\rho_0\big)
\cong
\big(\mathcal{L},0\big).
\]
\end{theorem}

Section 3 describes Labourie’s construction of twisted combinatorial cohomology with coefficients in the associated flat vector bundle. I derive an explicit combinatorial formula for infinitesimal deformations of flat connections on a closed surface \(\Sigma_g\). Using a good cover adapted to a filling ribbon graph \(\Gamma_g\), I describe the flat bundle through local flat trivializations and transition functions encoding its holonomy. This description admits the first-order variation of the associated discrete connection in terms of twisted combinatorial cochains. The resulting formula provides a concrete link between the differential-geometric description of infinitesimal deformations and the combinatorial model for twisted cohomology.\\\\
In (\cite{BillSymplecticNature}), Goldman constructed a natural symplectic form on the irreducible stratum of the character variety using the cup product in twisted cohomology and an invariant bilinear form on the Lie algebra. This motivates our study of the presymplectic structure on $\widetilde{\mathcal Y}_0$. Indeed, the space of noncentral $\mathrm{U}(1)$ representations is symplectic with the symplectic form $\omega^{\mathcal{Y}_0}$. The reducible noncentral stratum $\widetilde{\mathcal{Y}}_0$ is presymplectic under the Narasimhan-Atiyah-Bott-Goldman presymplectic form $\omega_{\mathsf{NABG}}$. The normal direction of $\mathcal{Y}_0$ in $\mathcal{R}$ is flat presymplecitic vector bundle $(\mathcal{N},\omega^{N})$ over the symplectic base $(\mathcal{Y}_0,\omega^{\mathcal{Y}_0})$. Embedding 
\[\big(\mathcal{Y}_0, \omega^{\mathcal{Y}_0}\big)\ \hookrightarrow\ \big(\mathcal{R},\omega_{\mathsf{NABG}}\big)\]
is a presymplectic embedding. 
\begin{theorem}{(\ref{PresymplecticEmbedding})}
\textit{
The presymplectic vector bundle $\big(\mathcal{N}, \omega^{N}\big)$ over the symplectic base $\big(\mathcal{Y}_0,\omega^{\mathcal{Y}_0}\big)$ determines a pre-symplectic embedding 
\[i:\ \big(\mathcal{Y}_0,\omega^{\mathcal{Y}_0}\big)\ \hookrightarrow\ \Big(\mathsf{Hom}\big(\pi_1(\Sigma_g),\mathrm{SU}(2)\big),\omega_{\mathsf{NABG}}\Big)\]
so that $i^*\omega_{\mathsf{NABG}}\ =\ \omega^{\mathcal{Y}_0}$.}
\end{theorem}
In particular, the Cup-Lie bracket capping with the fundamental of 2-cocycles in $H_2\big(\pi_1(\Sigma_g),\mathbb{R}\big)$ determines a splitting of the Zariski tangent space of the reducible noncentral in $\mathrm{Hom}\big(\pi_1(\Sigma_g),\mathrm{SU}(2)\big)$. Furthermore, the splitting is symmetric and nondegenerate on the vector bundle $\mathcal{N}$, and determinesns a subbundle of $\mathcal{N}$ over $\mathcal{Y}_0$ whose fiber is $\mathcal{C}_{\rho_0}$.
\begin{theorem}{(\ref{DecompositionOFConeIntoConeInTheNormalBundle})}
    \textit{Let $\rho_0$ be in $\mathcal{Y}_0$, the symmetric pairing $g_{\rho_0}$ determines a decomposition of the tangent quadratic cone $\mathcal{L}_{\rho_0}$ into the direct sum of a Euclidean space $\mathbb{R}^{2g}$ with a quadratic cone,
    \[\mathcal{L}_{\rho_0}\ \cong\ Z^1\big(\pi_1(\Sigma_g),\mathfrak{h}\big)\oplus \mathcal{C}_{\rho_0},\]
    where 
   \begin{equation}
       \label{QuadraicConeInNormalBundle}
       \mathcal{C}_{\rho_0}\ =\ \{u\in Z^1\big(\pi_1(\Sigma_g),\mathfrak{h}^{\perp}_{\mathrm{Ad}_{\rho_0}}\big):\ u\in\mathcal{L}_{\rho_0}\}.
   \end{equation}
    In particular, in the vector bundle $\mathcal{N}$, the quadratic cone $\mathcal{C}_{\rho_0}$ is a cone in the fiber $\mathcal{N}_{\rho_0}$.}
\end{theorem}

\begin{theorem}{(\ref{QuadraticConeFiberBundle})}
    \textit{The restriction of $g$ on $\mathcal{N}$ is symmetric nondegenerate. Therefore, the fiber bundle 
    \[\pi_C:\ \mathcal{C}\rightarrow\mathcal{Y}_0\]
    is a subbundle of $\mathcal{N}\rightarrow\mathcal{Y}_0$ so that $\pi_C^{-1}(\rho_0)=\mathcal{C}_{\rho_0}$.} 
\end{theorem}

In 1991, Sjamaar and Lerman developed the theory of stratified symplectic spaces (\cite{SjamaarStratified}), extending symplectic reduction to singular quotients and showing that the reduced space carries a natural stratification into symplectic manifolds. Let
\[
\chi
:=\mathrm{Hom}\bigl(\pi_1(\Sigma_g),\mathrm{SU}(2)\bigr)/\mathrm{SU}(2)
\]
be the character variety, with quotient map \(q:\mathcal R\to\chi\). For \(g\geq2\), the decomposition of \(\mathcal R\) descends to
\[
\chi
=
\chi_1\bigsqcup\chi_0\bigsqcup\chi^{\mathrm{irr}},
\]
where \(\chi_1=q(\mathcal Y_1)\), \(\chi_0=q(\widetilde{\mathcal Y}_0)\), and \(\chi^{\mathrm{irr}}=q(\mathcal R^{\mathrm{irr}})\) consist of central, noncentral reducible, and irreducible characters, respectively. Their closures satisfy
\[
\overline{\chi_1}=\chi_1,\qquad
\overline{\chi_0}=\chi_0\bigsqcup\chi_1,\qquad
\overline{\chi^{\mathrm{irr}}}=\chi.
\]
Taking connected components gives the corresponding strata. This construction describes the underlying orbit-type stratification without specifying a global algebra of smooth functions, a compatible Poisson bracket, or symplectic forms on the strata. By contrast, Sjamaar–Lerman’s singular reduction theorem (\cite{SjamaarStratified}) starts with a Hamiltonian group action and a moment map, and equips the reduced space with these compatible structures. The distinction therefore lies in the additional geometric structure established, rather than necessarily in the underlying partition into strata.\\\\
Jeffrey and Huebschmann developed finite-dimensional Hamiltonian
descriptions and stratified symplectic and Poisson structures for
moduli spaces of flat connections on surfaces
\cite{JeffreyExtendedModuliSpace,
      HuebschmannSingularetiesPossion,HuebschmannPoissonSU2}.
Building on these developments, our approach emphasizes the
geometry of the representation space along the noncentral
reducible stratum. I use a tubular neighborhood model in the normal bundle
$\mathcal{N}$ over $\mathcal{Y}_0$, making explicit the relationship
between the deformation cones, the stabilizer action, and the
local quotient geometry.\\\\
Combining the local normal form for moment maps developed by Sjamaar--Lerman (\cite{SjamaarStratified}) with the quasi-Hamiltonian double $\mathrm{SU}(2)\times\mathrm{SU}(2)$, I construct a model of an $\mathrm{SU}(2)$-orbit in $\mathsf{Hom}\big(\pi_1(\Sigma_g),\mathrm{SU}(2)\big)$.
\begin{proposition}{(\ref{SymplecticSlice})}
A neighborhood of the orbit $\mathcal{O}_{\rho_0}$ is equivariantly homeomorphic to a neighborhood of the zero section of
\[\mathcal{Z}_{\rho_0}\ =\ \mathrm{SU}(2)\times_{\mathrm{U}(1)}\Big(\mathrm{SU}(2)/\mathrm{U}(1)\times(\Psi_{\rho_0})_{dR}^{-1}(0)\Big),\]
 In particular, the embedding $\mathrm{SU}(2)/\mathrm{U}(1)\rightarrow \mathcal{Z}_{\rho_0}$ is isotropic and the normal bundle is 
\[\mathrm{SU}(2){\times}_{\mathrm{U}(1)}(\Psi_{\rho_0})_{dR}^{-1}(0).\]
\end{proposition}
Under the presymplectic embedding (\ref{PresymplecticEmbedding}), the cone bundle $\mathcal{C}\rightarrow\mathcal{Y}_0$ sits inside $\mathsf{Hom}\big(\pi_1(\Sigma_g),\mathrm{SU}(2)\big)$. The group $\mathrm{SU}(2)$ acts on the representation variety by conjugation. Meanwhile, $\mathrm{SU}(2)$ acts on the vector bundle $\mathcal{N}\rightarrow\mathcal{Y}_0$ in an equivariant way. 

\begin{theorem}{(\ref{Quotient Bundle Over The Kummer Variety})}
The vector bundle $\mathcal{N}\rightarrow\mathcal{Y}_0$ is $\mathrm{SU}(2)$-equivariant, and the quotient $\mathcal{N}/\mathrm{SU}(2)$ is a symplectic bundle over the noncentral classes in the Kummer variety $\mathcal{K}_0$. The fiber over $[\rho_0]$ is isomorphic to 
\[\Big(H^1\big(\pi_1(\Sigma_g),\mathfrak{h}_{\mathrm{Ad}_{\rho_0}}^{\perp}\big),\widetilde{\Omega}^N\Big).\]
\end{theorem}

Each fiber of the fiber bundle 
$\mathcal{N}/\mathrm{SU}(2)\rightarrow\mathcal{K}_0$ is a symplectic vector space, which is Hamiltonian under the $\mathrm{U}(1)$ action with the moment map in terms of the quadratic bracket.

\begin{theorem}{(\ref{A Reduced Symplectic Cone Bundle Over The Kummer Variety})}
    \textit{The reduced cone bundle 
    \[\widetilde{C}\ \rightarrow\ \mathcal{K}_0\ \]
    is a symplectic with the symplectic 2-form $\widetilde{\Omega}^{N}$, where the fiber at $[\rho_0]$ is the symplectic quadratic cone
    \[\big(\widetilde{\mathcal{C}}_{[\rho_0]}, \widetilde{\Omega}^N_{[\rho_0]}\big).\]
    The dimension of the reduced cone in the normal bundle at $[\rho_0]$ is $4g-6$.}
\end{theorem}

In 1998, Alekseev, Malkin, and Meinrenken introduced
quasi-Hamiltonian geometry~\cite{Alekseev1997LieGV}, establishing
a correspondence between compact quasi-Hamiltonian $G$-spaces
and Hamiltonian loop-group $LG$-spaces with proper moment maps.
A quasi-Hamiltonian space carries a group-valued moment map and
an invariant two-form satisfying a modified closedness condition
governed by the Cartan three-form, together with a minimal
degeneracy condition that ensures symplectic structures on
regular reductions. The noncentral representation classes in the Kummer variety $\mathcal{K}_0$ is the reducible noncentral stratum $\chi_0$, which is the q-Hamiltonian reduction of the Weyl group $W\cong\mathbb{Z}_2$. The Weyl group 
acts by identifying $\rho$ with its inverse $\rho^{-1}$.

\begin{proposition}{(\ref{Compatibility})}
    \textit{The q-Hamiltonian space $\big(\mathcal{Y}_0, \mathrm{Ad}_{W}, \omega, \mu|_{\mathcal{Y}_0}\big)$ sits in the tubular neighborhood $\mathcal{U}_{\mathcal{R}}$.
     In particular, $\omega^{\mathcal{Y}_0}$ is reduced to $\mathcal{Y}_0/W$ under the q-Hamiltonian reduction, $q^*\widetilde{\Omega}_{[\rho_0]}=\omega^{\mathcal{Y}_0}$ where
     \[q:\ \mathcal{Y}_0\ \rightarrow\ \chi_0.\] }
\end{proposition}

Therefore, the symplectic form of $\chi_0$ is the sum of the symplectic form of $\omega^{\mathcal{K}_0}$ and the symplectic form of the reduced cone bundle $\widetilde{\Omega}^{N}$.

\begin{theorem}{(\ref{SymplecticEmbedding})}
    The symplectic cone bundle $\widetilde{\mathcal{C}}\rightarrow\mathcal{K}_0$ sits in the character variety $\chi\big(\pi_1(\Sigma_g),\mathrm{SU}(2)\big)$, where $\mathcal{K}_0=\chi_0$. Therefore, the symplectic structure of $\chi_0$ is
    \[\widetilde{\Omega}_{\chi_0}\ =\ \omega^{\mathcal{K}_0}+\widetilde{\Omega}^N.\]
    \end{theorem}

\newpage
\section{Local Structure Of $\mathrm{Hom}
(\pi_1(\Sigma_g),\mathrm{SU}(2))$ At Reducible Representations}
Let $\Sigma_g$ be a closed surface of genus $g$, and
\[\mathrm{SU}(2)=\{\begin{bmatrix}
   a&b\\ -\bar{b}&\bar{a} 
\end{bmatrix}: |a|^2+|b|^2=1,\ a,\ b\in\mathbb{C}.\}\]
Reducible noncentral representations in $\mathsf{Hom}\big(\pi_1(\Sigma_g),\mathrm{SU}(2)\big)$ are
\[\begin{aligned}
     \pi_1(\Sigma_g)\ &\rightarrow\ \mathrm{SU}(2)\\
    x\ &\rightarrow\ \begin{bmatrix}
        \exp{\big(i\theta(x)\big)}&0\\0&\exp{\big(-i\theta(x)\big)} 
    \end{bmatrix}
\end{aligned}\]
where $\theta:\pi_1(\Sigma_g)\rightarrow\mathbb{R}$ is a group homomorphism and $\theta\not\equiv \mathbb{Z}\pi$. Let $i$ be the embedding 
\[\mathsf{Hom}\big(\pi_1(\Sigma_g),\mathrm{U}(1)\big)\ \hookrightarrow\ \mathsf{Hom}\big(\pi_1(\Sigma_g),\mathrm{SU}(2)\big).\]
Identify $\mathcal{Y}$ with the image of the embedding and  $\mathcal{Y}_0$ with the image of noncentral representations into $\mathrm{U}(1)$. Then, the closure of $\mathcal{Y}_0$ is $\mathcal{Y}$. In particular, the dimensions of the above spaces are
\[\begin{aligned}
    \dim_{\mathbb{R}}\mathrm{Hom}\big(\pi_1(\Sigma_g),\mathrm{SU}(2)\big)\ &=\ 6g-3.\\
    \dim_{\mathbb{R}}\mathcal{Y}\ =\ \dim_{\mathbb{R}}\mathcal{Y}_0\ &=\ 2g.
\end{aligned}\]
In this section, we introduce the algebraic tangent quadratic cone associated with a representation and explain how it arises from the local deformation theory of the representation variety. Starting from the space of infinitesimal deformations, described by group cohomology, we identify the quadratic obstruction that determines which first-order deformations extend to higher order. This construction provides a local model for the representation variety near a singular point and will serve as the basic framework for our later analysis of reducible representations.

\subsection{The formal deformation of a representation}
Before diving into the singularities, let's first review one of the main ideas of Goldman--Millson's work in \cite{BillMillsonDeforamtion}. First, we discuss the deformation theory of the $\mathbb{C}$-algebraic group $\mathrm{SL}(2,\mathbb{C})$, and then we apply it to $\mathrm{SU}(2)$.\\\\
Let $A_n=\mathbb{C}[[t]]/(t^{n+1})$ be an Artin local $\mathbb{C}$-algebra. The algebraic group $\mathrm{SL}(2, A_n)$ is a Lie group where the group structure is as follows, 
\[\mathrm{SL}(2,A_n)^0\ \rightarrow\mathrm{SL}(2,A_n)\ \rightarrow \mathrm{SL}(2,A_{n-1}).\]
$\mathrm{SL}(2,A_n)^0$ is the unipotent normal group of $\mathrm{SL}(2,A_n)$ so that
\begin{equation}
    \label{DeformationGroup}
    \mathrm{SL}(2,A_n)^0\ =\ \exp{\big(\mathfrak{sl}_2\otimes \mathfrak{m}_n\big)},
\end{equation}
and $\mathfrak{m}_n=(t^n)$ is the kernel of $A_{n}\rightarrow A_{n-1}$.\\\\
The set of formal deformations of $\rho_0$ is a compatible sequence $\{\rho_n\}_n$, so that
\[  \rho_n:\pi_1(\Sigma_g)\rightarrow \mathrm{SL}(2,A_n),\  q\circ\rho_n\ =\ \rho_{n-1},\ \forall\ n\geq1. \]
Explicitly, for aribitrary $x\in\pi_1(\Sigma_g)$,
\[\rho_n(x)\ =\ \exp{\big(u_n(x)\otimes(t^n)}\big)\cdot\rho_{n-1}(x),\]
where $u_n:\pi_1(\Sigma_g)\rightarrow\mathfrak{sl}_2$ is the map so that $\rho_n$ is a group homomorphism. The Lie bracket on $\mathfrak{g}\otimes \mathfrak{m}_n$ is defined by 
\[[\mathfrak{sl}_2\otimes \mathfrak{m}_n,\mathfrak{sl}_2\otimes\mathfrak{m}_n]\ =\ [\mathfrak{sl}_2,\mathfrak{sl}_2]\otimes\mathfrak{m}_n.\]
Let 
\[\mathcal{R}\big(\pi_1(\Sigma_g), \mathrm{SL}(2,\mathbb{C})\big)\]
denote the universal algebra associated with $\pi_1(\Sigma_g)$ and the algebraic group $\mathrm{SL}(2,\mathbb{C})$. The  representation scheme is 
\[\mathfrak{R}(\Sigma_g)\ :=\  \mathrm{spec}\Big(\mathcal{R}\big(\pi_1(\Sigma_g), \mathrm{SL}(2,\mathbb{C})\big)\Big),\]
and its closed points are naturally identified with the representation variety
\[\mathfrak{R}(\Sigma_g)(\mathbb{C})\ :=\ \mathrm{Hom}\big(\pi_1(\Sigma_g),\mathrm{SL}(2,\mathbb{C})\big)\]
 (\cite{SikoraCharacterVarieties}). For, $\rho\in \mathfrak{R}(\Sigma_g)(\mathbb{C})$, let $\mathfrak{m}_{\rho}$ be the corresponding maximal ideal and let $\hat{\mathcal{R}}_{\rho}$ be the completed local ring of the representation scheme at $\rho$. \\\\
 Let 
 \[A_1\ =\ \mathbb{C}[t]/(t^2)\]
 be the algebra of dual numbers. 
 Every analytic germ admits a canonical embedding into the Zariski tangent space $\mathrm{Hom}\big(\hat{\mathcal{R}}_{\rho}, A_1\big)$. Let $n=\dim_{\mathbb{C}}T_{\rho}\mathfrak{R}(\Sigma_g)$, then 
\[\mathfrak{m}_{\rho}/\mathfrak{m}_{\rho}^2\ \cong\ \mathbb{C}^n.\]
Then, the analytic germ of the $\mathrm{SL}(2,\mathbb{C})$-representation space at $\rho$ is defined by an ideal $\mathcal{J}_{\rho}$, so that
\[\hat{\mathcal{R}}_{\rho}\ =\ \mathbb{C}[[x_1,..., x_n]]/\mathcal{J}_{\rho}.\]
Suppose now that $\rho_0$ is a reducible representation, whose image is a compact subgroup of $\mathrm{SL}(2,\mathbb{C})$. By the Goldman-Millson quadraticity theorem, the analytic germ of the representation variety at $\rho_0$ is analytically equivalent to the germ of a quadratic cone. Equivalently, $\mathcal{J}_{\rho_0}$ is a homogeneous ideal which contains no homogeneous polynomials of degree 0 or 1.\\\\
Identify the Zariski tangent space of $\mathrm{Hom}\big(\pi_1(\Sigma_g),\mathrm{SU}(2)\big)$ at $\rho$ with the set 
\[\mathrm{Hom}\big(\hat{\mathcal{R}}_{\rho}, \mathbb{C}[t]/(t^2)\big).\] 
The tangent cone is the set of elements that are tangent to analytic paths starting at $\rho$ (\cite{Artin1968}). Hence, the tangent cone at $\rho$ is the image of the map 
\[\mathrm{Hom}\big(\hat{\mathcal{R}}_{\rho},\mathbb{C}[[t]]\big)\ \rightarrow\ \mathrm{Hom}\big(\hat{\mathcal{R}}_{\rho},\mathbb{C}[[t]]/(t^2)\big).\] 
Furthermore, the tangent quadratic cone equals the image of the map
\[\mathrm{Hom}\big(\hat{\mathcal{R}}_{\rho},\mathbb{C}[[t]]/(t^3)\big)\ \rightarrow\ \mathrm{Hom}\big(\hat{\mathcal{R}}_{\rho},\mathbb{C}[[t]]/(t^2)\big).\]
Hence, when $\rho_0$ is reducible whose image is a compact subgroup of $\mathrm{SL}(2,\mathbb{C})$, the deformation of $\rho_0$ is
\[\rho(x)\ =\ \exp{\big(u_1(x)\otimes t+u_2(x)\otimes t^2+\mathcal{O}(t^3)\big)}\cdot \rho_0(x),\ \forall\ x\in\pi_1(\Sigma_g).\]
where $u_i:\pi_1(\Sigma_g)\rightarrow\mathfrak{sl}_2 $ are the maps so that $\rho$ is a group homomorphism.\\\\ 
Now, consider the $\mathrm{SU}(2)$-representations. Recall that 
\[\mathfrak{sl}_2(\mathbb{C})\ \cong\ \mathfrak{su}(2)\otimes_{\mathbb{R}}\mathbb{C},\]
$\mathfrak{sl}_2(\mathbb{C})$ is the complexification of the real Lie algebra $\mathfrak{su}(2)$. Consider the conjugate-linear involution,
\[\begin{aligned}
    \kappa:\ \mathfrak{sl}_2(\mathbb{C})\ &\rightarrow\ \mathfrak{sl}_2(\mathbb{C})\\
    X\ &\mapsto\ -X^*=-\bar{X}^T. 
\end{aligned}\]
The real structure of $\mathfrak{sl}_2(\mathbb{C})$ is the Lie subalgebra that is fixed by $\kappa$,
\[\mathfrak{sl}_2(\mathbb{C})^{\kappa}\ =\ \{X\in\mathfrak{sl}_2(\mathbb{C}):\ \kappa(X)\ =\ X\}\ =\ \mathfrak{su}(2).\]
The corresponding conjugate-linear involution of the group $\mathrm{SL}(2,\mathbb{C})$ is
\[\begin{aligned}
    \widetilde{\kappa}:\ \mathrm{SL}(2,\mathbb{C})\ &\rightarrow\ \mathrm{SL}(2,\mathbb{C})\\
    A\ &\mapsto\ (\bar{A}^{T})^{-1}=(A^*)^{-1},
\end{aligned}\]
so that $\widetilde{\kappa}^2=\mathrm{Id}$. Moreover, 
\[\mathrm{SL}(2,\mathbb{C})^{\widetilde{\kappa}}\ =\ \{A\in\mathrm{SL}(2,\mathbb{C}):\ \widetilde{\kappa}(A)\ =\ A\}\ =\ \mathrm{SU}(2).\]
Let $I_n=(t^{n-1})/(t^n)$, $B_n=\mathbb{R}[[t]]/(t^{n+1})$ and $B_{n-1}=\mathbb{R}[[t]]/(t^n)$, then we obtain the following commutative diagram,
\[\begin{tikzcd}
	1 & {\exp{\big(\mathfrak{su}(2)\otimes_{\mathbb{R}}I_{n}^{\mathbb{R}}\big)}} & {\mathrm{SU}\big(2,B_n\big)} & {\mathrm{SU}\big(2, B_{n-1}\big)} & 1 \\
	1 & {\exp{\big(\mathfrak{sl}_2\otimes_{\mathbb{C}}I_{n}^{\mathbb{C}}\big)}} & {\mathrm{SL}\big(2,A_n\big)} & {\mathrm{SL}\big(2,A_{n-1}\big)} & 1
	\arrow[from=1-1, to=1-2]
	\arrow[from=1-2, to=1-3]
	\arrow[hook, from=1-2, to=2-2]
	\arrow[from=1-3, to=1-4]
	\arrow[hook, from=1-3, to=2-3]
	\arrow[from=1-4, to=1-5]
	\arrow[hook, from=1-4, to=2-4]
	\arrow[from=2-1, to=2-2]
	\arrow[from=2-2, to=2-3]
	\arrow[from=2-3, to=2-4]
	\arrow[from=2-4, to=2-5]
\end{tikzcd}\]
The columns are the inclusions of the fixed-point subgroups and fixed-point Lie algebras into their complexifications. From the injective form of the short five lemma, we obtain the following lemma.
\begin{lemma}
    The group homomorphism
    \[\mathrm{SU}(2,B_n)\ \rightarrow\ \mathrm{SL}(2,A_n)\]
    is injective.
\end{lemma}
Let $\rho:\pi_1(\Sigma_g)\rightarrow\mathrm{SU}(2)$ be a representation, and let $\mathfrak{su}(2)_{\mathrm{Ad}_{\rho}}$ be the corresponding $\pi_1(\Sigma_g)$-module. The Zariski tangent space of $\mathrm{Hom}\big(\pi_1(\Sigma_g),\mathrm{SU}(2)\big)$ at $\rho$ is identified with
\[ Z^1\big(\pi_1(\Sigma_g),\mathfrak{su}(2)_{\mathrm{Ad}_{\rho}}\big).\]
For $u\in Z^1\big(\pi_1(\Sigma_g),\mathfrak{su}(2)_{\mathrm{Ad}_{\rho}}\big)$, define a 2-cocycle $[u,u]$ so that for all $x,y\in\pi_1(\Sigma_g)$,
\begin{equation}
    \label{DefinitionOfBracket}
    [u,u](x,y)\ :=\ [u(x),\mathrm{Ad}_{\rho(x)}u(y)]_{\mathfrak{su}(2)}.
\end{equation}
Indeed, the cohomology class of $[u,u]$ is the primary obstruction to extend the infinitesimal deformation $u$ to the second order. Denote the algebraic tangent quadratic cone at a $\rho$ to be $\mathcal{L}_{\rho}$, which is the space of 1-cocycles $u$ so that the class of $[u,u]$ equals to zero in $H^2(\pi_1(\Sigma_g),\mathfrak{su}(2)_{\mathrm{Ad}_{\rho}})$ (\cite{BillMillsonDeforamtion}). Let $\delta$ be the co-boundary map of the $\pi_1(\Sigma_g)$-complex $C^{\bullet}\big(\pi_1(\Sigma_g), \mathfrak{su}(2)_{\mathrm{Ad}_{\rho}}\big)$. For any $u\in\mathcal{L}_{\rho}$, there exists a 1-cochain 
\[v:\ \pi_1(\Sigma_g)\ \rightarrow\ \mathfrak{su}(2)\]
satisfying 
\begin{equation}
    \label{TangentQuadraticConeEq}
    [u,u]\ =\ -\frac{1}{2}\delta v,
\end{equation}
where for any $x,y\in\pi_1(\Sigma_g)$,
\[\delta v(x,y)\ =\ \mathrm{Ad}_{\rho(x)}v(y)-v(xy)+v(x).\]

\subsection{The tangent quadratic cone at a reducible noncentral representation}
Let $\rho_0\in\mathcal{Y}_0$, then 
\[\rho_0:\ \pi_1(\Sigma_g)\ \rightarrow\ \mathrm{U}(1)\ \xrightarrow{i}\ \mathrm{SU}(2).\] 
Associated with $\rho_0$ is the flat principal $\mathrm{SU}(2)$-bundle 
\[P_{\rho_0}\ :=\ \widetilde{\Sigma}_g\times_{\rho_0}\mathrm{SU}(2)\ \rightarrow\ \Sigma_g.\]
The holonomy of its flat connection is contained in $\mathrm{U}(1)$, the structure group of $P_{\rho_0}$ admits a reduction from $\mathrm{SU}(2)$ to $\mathrm{U}(1)$.
Suppose $x_0$ is a base point of $\Sigma_g$, and $\widetilde{x}_0$ lies over $x_0$ in the universal cover $\widetilde{\Sigma}_g$ of $\Sigma_g$. $\pi_1(\Sigma_g, x_0)$ is a group of deck transformations on
\[\big(\widetilde{\Sigma}_g, \widetilde{x}_0\big)\ \rightarrow\ \big(\Sigma_g,x_0\big).\]
For any $\gamma\in\pi_1(\Sigma_g, x_0)$, $\gamma$ acts on the product bundle $\widetilde{\Sigma}_g\times\mathrm{SU}(2)\rightarrow\widetilde{\Sigma}$ by
\[(\tilde{x}_0, g)\ \mapsto\ (\gamma\cdot\tilde{x}_0,\ \rho(\gamma)\cdot g).\]
 The associated bundle 
\begin{equation}
\label{AssociatedHomogeneousBundle}
    E_{\rho_0}\ :=\ P_{\rho_0}\times_{\mathrm{U}(1)}\big(\mathrm{SU}(2)/\mathrm{U}(1)\big)\ \rightarrow\ \Sigma_g
\end{equation}
with the structure group $\mathrm{U}(1)$ admits a parallel section with respect to the flat connection (\cite{KobayashiDifferentialGeometry1}),
\begin{equation}
    \label{ParallelSection}
    \sigma_0:\ \Sigma_g\ \rightarrow\  E_{\rho_0}.
\end{equation}
Suppose $\mathcal{B}$ is the Killing form of $\mathfrak{su}(2)$, so that as a vector space
\[\mathfrak{su}(2)\ =\ \mathfrak{h}+\mathfrak{h}^{\perp},\]
and $[\mathfrak{h},\ \mathfrak{h}]=0,\ [\mathfrak{h},\mathfrak{h}^{\perp}]\subset\mathfrak{h}^{\perp}$. Then $\mathfrak{h}^{\perp}$ is a representation vector space of $\mathrm{U}(1)$, so that
\[\mathrm{ad}E_{\rho_0}=P_{\rho_0}\times_{\mathrm{U}(1)}\mathfrak{h}^{\perp}\ \rightarrow\ \Sigma_g\]
is an associated vector bundle with the structure group $\mathrm{U}(1)$.  
\[Z^1\big(\pi_1(\Sigma_g),\mathfrak{h}\big)\ \text{and}\  Z^1\big(\pi_1(\Sigma_g),\mathfrak{h}^{\perp}_{\mathrm{Ad}_{\rho_0}}\big)\]
are two subspaces of the Zariski tangent space
\begin{equation}
    \label{ZariskiTangentSpaceOfRepresentationVariiety}
    \mathcal{T}_{\rho_0}\ :=\ T_{\rho_0}\mathrm{Hom}\big(\pi_1(\Sigma_g),\mathrm{SU}(2)\big)
\end{equation}
\begin{enumerate}
\item If $u\in Z^1\big(\pi_1(\Sigma_g),\ \mathfrak{h}\big)$, where $\mathfrak{h}$ is abelian, the 2-cocycle
\[[u,u]\ \equiv\ 0.\] 
Taking $v\equiv0$, any $u$ in $Z^1\big(\pi_1(\Sigma_g), \mathfrak{h}\big)$ satisfies the equation (\ref{TangentQuadraticConeEq}). Hence, 
\[Z^1\big(\pi_1(\Sigma_g),\mathfrak{h}\big)\ \subset\ \mathcal{L}_{\rho_0}.\]
Moreover, $\dim_{\mathbb{R}}Z^1\big(\pi_1(\Sigma_g), \mathfrak{h}\big)=2g$.
\item If $u$ is in $Z^1\big(\pi_1(\Sigma_g), \mathfrak{h}^{\perp}_{\mathrm{Ad}_{\rho_0}}\big)$, then $[u,\ u]$ is a 2-cocycle taking values in $\mathfrak{h}$ (\ref{TangentQuadraticConeEq}), and for any $x,y\in\pi_1(\Sigma_g)$,
\[[u(x), \mathrm{Ad}_{\rho_0(x)}u(y)]\ =\ -2\big(v(x)-v(y)+v(xy)\big).\]
Construct a 1-cochain $v\in C^1\big(\pi_1(\Sigma_g),\mathfrak{h}\big)$, by setting 
\[v(a_i)\ =\ v(b_i)\ =\ 0,\ \forall\ i=1,2,...,g, \]
 and determining its values on the remaining group elements recursively using the above equation. Hence, the dimension of
 \[\mathcal{C}_{\rho_0}\ :=\ Z^1\big(\pi_1(\Sigma_g), \mathfrak{h}^{\perp}_{\mathrm{Ad}_{\rho_0}}\big)\ \cap \mathcal{L}_{\rho_0}\]
 is $4g-3$.
\end{enumerate}

\newpage
\section{Tangent Quadratic Cone in Combinatorial Version}
The algebraic tangent quadratic cone at reducible but noncentral representations is solved in terms of the deformation of the flat principal bundles with respect to the corresponding flat connections. In particular, the deformed principal bundle with the Artin ring-valued structure group provides a geometric view of the representation variety. \\\\
Restricting a flat vector bundle to a good cover adapted to a filling ribbon graph \(\Gamma_g\) describes the bundle through local flat trivializations and locally constant transition functions encoding its holonomy. This local description identifies the twisted de Rham cohomology \(H^\bullet_{\mathrm{dR}}(\Sigma_g;\mathrm{ad}P)\) with the corresponding combinatorial twisted cohomology $H^{\bullet}\big(\Gamma_g,\mathrm{ad}P\big)$. The first cohomology describes infinitesimal deformations of the associated discrete connection modulo infinitesimal gauge transformations. The full local deformation problem is governed by the Maurer–Cartan equation
\[
d_Aa+\frac12[a\wedge a]=0,
\qquad a\in\Omega^1(\Sigma_g;\operatorname{ad}P),
\]
with solutions considered modulo gauge transformations.

 \subsection{The combinatorial twisted cohomology theory}
Inspired by Labourie's work (\cite{FranciosRepresentation}), we fix a filling Ribbon graph $\Gamma_g$ to be a single-vertex petal shape with $2g$ edges. The graph is oriented, and $\Gamma_g^0$, $\Gamma_g^1$ and $\Gamma_g^2$ are sets of oriented vertices, oriented edges and oriented faces, respectively. In particular, $\Gamma_g^1=\{v_0\}$, $\Gamma_g^2=\{f\}$, and
\[\Gamma_g^1=\{a_1,b_1,...,a_g,b_g\}.\]
Any element $\alpha\in \Gamma_g^i$, let $\overline{\alpha}$ be the element $\alpha$ with the opposite direction.  The boundaries are:
\[\begin{aligned}
    \partial v_0\ &=\ 0; \\
    \partial l\ &=\ \big(l_+,\overline{l_-}\big),\ \forall\ l\in \Gamma_g^1\\
    \partial f\ &=\ \big(a_1,b_1,\overline{a_1},\overline{b_1},...,a_g,b_g,\overline{a_g},\overline{b_g}\big) 
\end{aligned}\]
For any $\alpha\in \Gamma_g$, assign $\alpha$ an contractible open set $U_{\alpha}\subset\Sigma_g$ containing the interior $\hat{\alpha}$ of $\alpha$. Let
\[W_1\ :=\ \bigsqcup_{i=1}^gU_{a_i}\bigsqcup_{i=1}^gU_{b_i}.\]
For each pair $(\alpha, \beta)$ so that $\alpha\in\partial\beta$, take open sets $U_{\alpha, \beta}$ that are homeomorphic to discs. Each $U_{\alpha}$ satisfies the following conditions.
\begin{itemize}
    \item $U_{\overline{\alpha}}= U_{\alpha}$ and $U_{\overline{l}, \overline{f}}=U_{l,f}$.
    \item $\forall \alpha, \beta\in \Gamma_g^1$ so that $\alpha\neq\beta$ or $\alpha\neq\bar{\beta}$, 
    \[\overline{U}_{\alpha}\cap \overline{U}_{\beta}\ =\ \emptyset;\]
    \item For every edge $l$, 
    \[U_l\cap U_{v_0}=U_{l_+, l}\bigsqcup U_{\overline{l_-},l};\]
    \item For the single face $f$, 
    \[W_1\cap U_{f}\ =\ \bigsqcup_{i=1}^gU_{a_i,f}\bigsqcup_{i=1}^gU_{b_i,f}.\]
\end{itemize}
Let $\mathcal{U}^{\Gamma_g}$ be the collection of all such open sets relative to $\Gamma_g$, then $\mathcal{U}^{\Gamma_g}$ is a good covering of $\Sigma_g$. \\\\
Now, the flat connection $\triangledown_{\rho_0}$ has a combinatorial description relative to $\Gamma_g$, which is still denoted by $\triangledown_{\rho_0}$. Assigning the group $\mathrm{SU}(2)$ to each edge $l$ , the combinatorial connection is a locally constant map
\[\begin{aligned}
    \triangledown_{\rho_0} :\ \Gamma_g^1&\ \rightarrow\ \mathrm{SU}(2)\\
    l&\ \mapsto g_l,\ g_{\bar{l}}\ =\ g_l^{-1}.
\end{aligned}\]
Indeed, for an edge $l\in E$, $l_-=l_+=v_0$, the closure of $l$ is
\[\check{l}\ =\ l\cup \{v_0\},\]
so that $\check{l}$ is a representative of a homotopy class in $\pi_1(\Sigma_g,v_0)$. Set
\[g_{l}\ =\ \rho_0(\check{l}),\]
so that $g_{\overline{l}}\ =\ \rho_0(\check{l})^{-1}$. Assign $\mathrm{SU}(2)$ to the single vertex $v_0$. Meanwhile, the group $\mathrm{SU}(2)^{\{v_0\}}=\mathrm{SU}(2)$ is the gauge group of $\Gamma_g$, and acts on $\mathrm{SU}(2)^{\Gamma_g^1}$ via
\[\{h_{v_0}\}\cdot \{g_l\}_{l\in \Gamma_g^1}:=\ \{h_{v_0}g_lh_{v_0}^{-1}\}_{l\in \Gamma_g^1}.\]
The bundle map of the associated vector bundle $\mathrm{ad}P$ is $\pi_1: \mathrm{ad}P\rightarrow\Sigma_g$, where
\[\mathrm{ad}P\ :=\ P\times_{\mathrm{SU}(2)}\mathfrak{su}(2).\]
Now, we define the complex, 
\[0\ \rightarrow\ C_{\Gamma_g}^0(\mathrm{ad}P_{\rho_0})\ \xrightarrow{\delta^0}\ C_{\Gamma_g}^1(\mathrm{ad}P_{\rho_0})\ \xrightarrow{\delta^1}\ C_{\Gamma_g}^2(\mathrm{ad}P_{\rho_0})\ \rightarrow\ 0\]
where
\[C^i_{\Gamma_g}(\mathrm{ad}P_{\rho_0})\ :=\ \{c^i:\ \Gamma_g^i\rightarrow\bigsqcup_{\alpha\in \Gamma_g^i}(\mathrm{ad}P_{\rho_0}|_{U_{\alpha}}): c^i(\alpha)\in\mathrm{ad}P_{\rho_0}|_{U_{\alpha}},\ c^i(\bar{\alpha})\ =\ -c^i(\alpha)\}.\]
The coboundary is 
\[\begin{aligned}
    \delta^i:\ C_{\Gamma_g}^i&\ \rightarrow\ C_{\Gamma_g}^{i+1}\\
    c^{i}&\ \mapsto\ (\delta^ic^i)(\beta_{i+1}):=\sum_{\alpha_i\in\partial\beta_{i+1}}i_{\alpha_i,\beta_{i+1}}c^i(\alpha_i).
\end{aligned}\]
where $i_{\alpha,\beta}$ is the assignment of parallel sections along $U_{\alpha}$ to the unique parallel section along $U_{\beta}$ coincides with $u$ on $U_{\alpha,\beta}$.\\\\
The open set $U_{\alpha}$ is contractible, the locally constant sections of $\mathrm{ad}P_{\rho_0}$ over $U_{\alpha}$ form a vector space. In particular, at the single vertex $v_0$, 
\[V_{v_0}\ =\ \{s:U_{v_0}\rightarrow\mathrm{ad}P_{\rho_0}:\ \nabla_{\rho_0}s=0\}.\]
A canonical local trivialization of $\mathrm{ad}P_{\rho_0}$ at $U_{v_0}$ is
\[\begin{aligned}
    \tau:\ U_{v_0}\times V_{v_0}\ &\rightarrow\ \mathrm{ad}P_{\rho_0}|_{U_{v_0}}\\
    (x, s)\ &\mapsto\ s(x)
\end{aligned}.\] 
Hence, at point $x\in U_{v_0}$ $\tau_x=\mathrm{ev}_x$, the evaluation map at $x$.
 \begin{enumerate}
     \item 
     If $c\in C_{\Gamma_g}^0(\mathrm{ad}P_{\rho_0})$, so that $c(v_0)=\xi$, for some $\xi\in\mathfrak{su}(2)$. Then
     \[\begin{aligned}
         (\delta^0c)(l)\ &=\ i_{l^+, l}(c(l_+))+i_{\overline{l_-}, l}(c(\overline{l_-}))\\
         &=\ i_{l^+, l}(c(l_+))-i_{\overline{l_-}, l}(c(l_-)),
     \end{aligned}\]
where $l$ is an edge, and $l_+=l_-=v_0$. Identify the fame at $l^+$ with the local trivialization at $v_0$,
\[i_{\overline{l_-},l}\ =\ \mathrm{Id},\ i_{l_+,l}\ =\ g_l\]
Hence, let $c(l_-)=\xi$, for some $\xi\in\mathfrak{su}(2)$,
\[ (\delta^0c)(l)\ =\ g_l\cdot c(l_+)- c(l_-)\ =\ \mathrm{Ad}_{g_l}\xi-\xi.\]
 \item If $c\in  C_{\Gamma_g}^1(\mathrm{ad}P_{\rho_0})$, then
   \begin{equation}
       \label{OneCocycleConditionCombinatorial}
       \begin{aligned}
           (\delta^1c)(f)\ &=\ \sum_{i=1}^g\big(i_{a_i,f}c(a_i)+i_{b_i,f}c(b_i)+i_{\overline{a_i},f}c(\overline{a_i})+i_{\overline{b_i},f}c(\overline{b_i})\big)\\
       &=\ \sum_{i=1}^g\big(i_{a_i,f}c(a_i)+i_{b_i,f}c(b_i)-i_{\overline{a_i},f}c(a_i)-i_{\overline{b_i},f}c(b_i)\big)\\
       &=\ \sum_{i=1}\big((\mathrm{Id}-g_{a_i}g_{b_i}g_{\overline{a}_i})c(a_i)+(g_{a_i}-g_{a_i}g_{b_i}g_{\overline{a}_i}g_{\overline{b_i}})c(b_i)\big)\\
       &=\ \sum_{i=1}\big((\mathrm{Id}-g_{b_i})\cdot c(a_i)+(g_{a_i}-\mathrm{Id})\cdot c(b_i)\big),
       \end{aligned}    
   \end{equation}
   where $\partial f=\prod_{i=1}^g[a_i,b_i]$, the product of commutators. Consequently, if $c$ is a 1-cochain, then $c$ will satisfy that
   \[(\delta^1c)(f)\ =\ c(\partial f)\ =\ 0.\]
 \end{enumerate}

\subsection{The tangent quadratic tangent cone}
Let $\triangledown_{\rho_0}$ be the flat connection corresponding to $\rho_0$ on the principal bundle $\pi: P\rightarrow\Sigma_g$, and
\[\triangledown_{\rho_0}\ =\ d+\omega_0,\]
where $\omega_0$ is the corresponding connection form. In general, as the  connection form of a $\mathrm{SU}(2)$-bundle, $\omega_0$ has a decomposition
\[\omega_0\ =\ (\omega_0)_{\mathfrak{h}}+(\omega_0)_{\mathfrak{h}^{\perp}}.\]
The representation $\rho_0$ is reducible noncentral in $\mathsf{Hom}\big(\pi_{\Sigma_g},\mathrm{SU}(2)\big)$, the flat connection $\triangledown_{\rho_0}$ determines a section  of the associated bundle (\ref{AssociatedHomogeneousBundle}) $\sigma_0:\Sigma_g\rightarrow E_{\rho_0}$. Identify $E_{\rho_0}$ with $P_{\rho_0}/\mathrm{U}(1)$, the bundle $P_{\rho_0}\rightarrow E_{\rho_0}$ has the structure group $\mathrm{U}(1)$. Indeed, the section $\sigma_0$ is a parallel section of the bundle 
\[\Sigma_g\ \rightarrow\ P_{\rho_0}/H.\]
Now, replace the representation $\rho_0$ with the section $\sigma_0$ of the principal $\mathrm{SU}(2)$-bundle $P$ with respect to the flat connection $\nabla_{\rho_0}$. For edges $a_i$ and $b_j$, $\check{a}_i$ and $\check{b}_j$ are representatives of generators of $\pi_1(\Sigma_g)$. Restricting the section $\sigma_0$ to $U_{l}$, a discrete connection on the filling Ribbon graph $\Gamma_g$ is a locally constant sections of $P_{\rho_0}|_{U_l}$'s, for all $l\in\Gamma_g^1$, 
\[\begin{aligned}
    \sigma_0:\ U_{l}\ &\rightarrow\ \pi^{-1}(U_{l})\cong U_{l}\times\mathrm{SU}(2)\\
    x\ &\mapsto\ \Big(x,p_2\big(\sigma_0(x)\big)\Big),
\end{aligned}\]
where $p_2$ is the projection on the fiber over $x$. Since $U_l$ is contractible, $\sigma_0$ is determined by its value at a point $x\in U_{l}$. Identify $p_2\circ\sigma_0(x)$ with $\sigma_0$,
then $\sigma_0(l)=\rho_0(\check{l})$.\\\\
Let $P_A$ is the associated bundle
\[P_A\ :=\ P\times_{\mathrm{SU}(2)}\mathrm{SU}(2, A_2),\]
where the bundle map is $\pi_A:\ P_{A}\ \rightarrow\ \Sigma_g$. Furthermore, $P_A\rightarrow P$ is a principal fiberation with the structure group $\mathrm{SU}(2,A_2)^0$ (\cite{BillMillsonDeforamtion}). Identify a locally constant section $\sigma_0$ with the discrete connection corresponding to $\nabla_{\rho_0}$ on $P$, the infinitesimal deformation space of $\sigma_0$ is isomorphic to the space $Z^1\big(\Gamma_g,\mathsf{ad}P_{\rho_0}\big)$. Suppose $C^{\bullet}_{\Gamma_g}(\mathfrak{h})$ and $C_{\Gamma
_g}^{\bullet}\big(\mathfrak{h}_{\mathrm{Ad}_{\rho_0}}^{\perp}\big)$ are combinatorial twisted complexes of the flat vector bundles $P\times_{\rho_0}\mathfrak{h}$ and $\mathrm{ad}E_{\rho_0}$ respectively. In particular, the complex $C^{\bullet}_{\Gamma_g}(\mathfrak{h})$ is
\[0\ \rightarrow\ \mathfrak{h}\ \xrightarrow{0}\ \mathfrak{h}^{2g}\ \xrightarrow{0}\ \mathfrak{h}\ \rightarrow\ 0.\]
Hence,
\[\begin{aligned}
    Z^1\big(\Gamma_g,\mathfrak{h}\big)\ &\cong\ \mathfrak{h}^{2g},\\ Z^2\big(\Gamma_g,\mathfrak{h}\big)\ &=\ \mathfrak{h}\ =\ Z^0\big(\Gamma_g,\mathfrak{h}\big),\\ B^1\big(\Gamma_g,\mathfrak{h}\big)\ &=\{0\}\ =\ B^2\big(\Gamma_g,\mathfrak{h}\big).
\end{aligned}\]
where $Z^1\big(\Gamma_g,\mathfrak{h}^{\perp}_{\mathrm{Ad}}\big)$ is the space of 1-cocycles, and $B^1\big(\Gamma_g,\mathfrak{h}^{\perp}_{\mathrm{Ad}_{\rho_0}}\big)$ is the space of 1-coboundaries.\\\\
Now, let $l$ be an edge of $\Gamma_g$, and  $U_l\in\mathcal{U}^{\Gamma_g}$ is the open subset containing the interior $\hat{l}$ of $l$. Over $U_l$, assume $\sigma_0=e\mathrm{U}(1)$, the identity coset in $\mathrm{SU}(2)/\mathrm{U}(1)$. Let $U_0$ be the open set of $\mathrm{SU}(2)/\mathrm{U}(1)$ where
\[\mathrm{SU}(2)\rightarrow\mathrm{SU}(2)/\mathrm{U}(1)\]
is an analytic diffeomorphism. Lift $\sigma_0$ to $\sigma_1$ that takes values in $\mathrm{SU}(2)$, 
\[\begin{aligned}
    \sigma_1: U_l\ &\rightarrow\ U_l\times \mathrm{SU}(2)\\
    x\ &\mapsto\ \big(x, \sigma_1(x)\big).
\end{aligned}\]
Indeed, the section $\sigma_1$ is a section of the bundle $P_{A}$, so that
\[\sigma_1(l)\ =\ \exp{\big(c_1(l)\otimes t+c_2(l)\otimes t^2+\mathcal{O}(t^3)\big)}\cdot\sigma_0(l),\]
where $c_1, c_2: \Gamma_g^1\rightarrow\mathfrak{su}(2)$ are locally constant sections of the associated vector bundle $\mathrm{ad}P_{\rho_0}$. From the fact that $c_1,c_2$ correspond to the deformation of the discrete connection $\sigma_0$, the maps $c_1, c_2$ are 1-cocycles.
Now, write $\sigma_1$ as a matrix in $\mathrm{SU}(2)$ over each edge,  
 \[\sigma_1(l)\ =\ \frac{1}{\sqrt{1+|z(t,l)|^2}}\begin{bmatrix}
    1&z(t,l)\\-\overline{z(t,l)}&1
\end{bmatrix}.\]
Then the connections form on $U_l$ is
    \[\begin{aligned}
         \sigma_1(l)^{-1}d\sigma_1(l)\ =
    \frac{1}{\sqrt{1+|z(l,t)|^2}}\begin{bmatrix}
        \frac{-\overline{z(l,t)}dz(l,t)+z(l,t)d\overline{z(l,t)}}{2}&dz(l,t)\\-d\overline{z(l,t)}&\frac{-z(l,t)d\overline{z(l,t)}+\overline{z(l,t)}dz(l,t)}{2}
    \end{bmatrix}.
\end{aligned}\]
In terms of the polar coordinates,
\[z(t,l)\ =\ \ r(t,l)\cos{\phi(t,l)}\ +\ ir(t,l)\sin{\phi(t,l)},\]
the deformation of the connection form is 
\[\begin{aligned}
    \sigma_1^{-1}d\sigma_1\ =&\ \big(\frac{-r(t,l)\sin{\phi(t,l)}\phi'(t,l)+r'(t,\gamma)\cos{\phi(t,l)}}{1+r(t,l)^2}X\\
    &+\frac{r(t,l)\cos{\phi(t,l)}\phi'(t,l)+r'(t,\gamma)\sin{\phi(t,l)}}{1+r(t,l)^2}Y-\frac{r(t,l)^2\phi'(t,l)}{1+r(t,l)^2}H\big)dt.
\end{aligned}\]
Let $Z^1\big(\Gamma_g,\mathfrak{h}\big)$ be the space of combinatorial 1-cocycles.
\begin{itemize}
    \item If both $\phi'=r'=0$, then $\sigma_1=\sigma_0$;
    \item If $r(t,l)^2\phi'(t,l)\neq0$, and 
    \[\begin{aligned}
        -r(t,l)\sin{\phi(t,l)}\phi'(t,l)+r'(t,\gamma)\cos{\phi(t,l)}\ &=\ 0\\
        r(t,l)\cos{\phi(t,l)}\phi'(t,l)+r'(t,\gamma)\sin{\phi(t,l)}\ &=\ 0.
    \end{aligned}\]
    Solve the equation, 
    \[r'\ =\ \phi\ =\ 0.\]
    In this case, 
    \[\sigma_1\ =\ \frac{1}{\sqrt{1+r^2}}\begin{bmatrix}
        1&r\\-r&1
    \end{bmatrix},\ r\in\mathbb{R}-\{0\},\]
    and the infinitesimal deformation space is $Z^1\big(\Gamma_g,\mathfrak{h}\big)$.
    \item  Otherwise, the coefficient of $H$ equals to zero,
    \[r(t,l)^2\phi'(t,l)=0.\]
    Equivalently, the angle $\phi$ is independent of $t$. Set $\phi=\phi_0$. In this case, 
\[\begin{aligned}
    \sigma_1\ =\ \frac{1}{\sqrt{1+r(t)^2}}\begin{bmatrix}
    1& r(t)\cos{\phi_0}+ir(t)\sin{\phi_0}\\
    -r(t)\cos{\phi_0}+ir(t)\sin{\phi_0}&1
\end{bmatrix}
\end{aligned}\]
and the infinitesimal deformation is 
\[\frac{r'}{1+r(t)^2}\big(\cos{\phi_0}X+\sin{\phi_0} Y\big),\]
which is in $Z^1\big(\Gamma_g,\mathfrak{h}_{\mathrm{Ad}_{\rho_0}}^{\perp}\big)$. Let $f$ be any smooth function in $0$ in the 1-cocycle space, and 
\[r(t)\ =\ \tan{\int_0^tf(s)}ds.\]
Then, $\frac{r'}{1+r^2}=f$. Hence, the deformation space consists of the 1-cocycles 
\[f(t)\big(\cos{\phi_0}X+\sin{\phi_0}Y\big).\]
\end{itemize}

\newpage
\section{The Presymplectic Reducible Noncentral Stratum In The Representation Variety}
The set of central representations is $\mathcal{Y}_1$, which consists of $2^{2g}$ isolated points in $\mathrm{Hom}\big(\pi_1(\Sigma_g),\mathrm{SU}(2)\big)$. The stratum $\widetilde{\mathcal{Y}}_0$ is the set of reducible noncentral representations in $\mathrm{Hom}\big(\pi_1(\Sigma_g),\mathrm{SU}(2)\big)$, which are representations conjugate to representations in $\mathcal{Y}_0$ under $\mathrm{SU}(2)$. The space $\mathrm{Hom}\big(\pi_1(\Sigma_g),\mathrm{SU}(2)\big)^{\mathrm{irr}}$ is the space of irreducible representations.
\begin{proposition}{\cite{SjamaarStratified}}
\label{StratifiedRepresentationVariety}
The representation space $\mathcal{R}$, viewed as a closed subset
of $\mathrm{SU}(2)^{2g}$ via the standard presentation of
$\pi_1(\Sigma_g)$, admits the orbit-type decomposition
\[
\mathcal{R}
=
\mathcal{Y}_1
\bigsqcup
\widetilde{\mathcal{Y}}_0
\bigsqcup
\mathcal{R}^{\mathrm{irr}}.
\]
The stabilizers under conjugation on these three pieces are
conjugate to $\mathrm{SU}(2)$, $\mathrm{U}(1)$, and
$\{\pm I\}$, respectively. Their closures in $\mathcal{R}$ satisfy
\[
\overline{\mathcal{Y}_1}=\mathcal{Y}_1,
\qquad
\overline{\widetilde{\mathcal{Y}}_0}
=
\widetilde{\mathcal{Y}}_0\bigsqcup\mathcal{Y}_1,
\qquad
\overline{\mathcal{R}^{\mathrm{irr}}}=\mathcal{R}.
\]
In particular, under the frontier order
$S\preceq T$ defined by $S\subseteq\overline{T}$, these pieces form
the chain
\[
\mathcal{Y}_1
\prec
\widetilde{\mathcal{Y}}_0
\prec
\mathcal{R}^{\mathrm{irr}}.
\]
\end{proposition}

\begin{remark}
\textit{The disjoint union of  $\mathsf{Hom}\big(\pi_1(\Sigma_g),\mathrm{SU}(2)\big)$ in the above proposition automatically satisfies the Whitney condition. The infinitesimal deformation space determines the analytic paths $\gamma^+$ starting at each point on the representation variety. For any representation $\rho$, the infinitesimal deformation is a 1-cocycle taking values in the corresponding twisted flat vector bundle. Hence, the path from $\rho$ is not reversed. If a representation $\rho$ connects to $\rho_0\in\mathcal{Y}$, then $\rho$ will be reducible and be in $\mathcal{Y}_0$.}
\end{remark}

\subsection{Incomplete vector fields on the reducible noncentral stratum}
This section discusses the differential geometry of the representation variety in terms of the group cohomology theory. Let $Z^1\big(\pi_1(\Sigma_g),\mathfrak{su}(2)_{\mathrm{Ad}}\big)$ be the vector bundle over $\mathcal{Y}_0$ in $\mathsf{Hom}\big(\pi_1(\Sigma_g),\mathrm{SU}(2)\big)$. Identify the representation variety $\mathsf{Hom}\big(\pi_1(\Sigma_g), \mathrm{U}(1)\big)$ with $\mathrm{U}(1)^{2g}$. Suppose the word map is
\[\begin{aligned}
   w:\ \mathrm{U}(1)^{2g}\ &\rightarrow\ \mathrm{U}(1)\\
    (x_1,y_1,...,x_g, y_g)\ &\mapsto\ x_1y_1x_1^{-1}y_1^{-1}\cdots x_gy_gx_g^{-1}y_g
^{-1}.
\end{aligned}\]
A representation $\rho_0$ is identified with a unique point $\big(\rho_0(a_1), \rho_0(b_1),\cdots,\rho_0(a_g),\rho_0(b_g)\big)$ in $\mathrm{U}(1)^{2g}$. Indeed, the word map $w$ is constant so that $w\equiv1$. For any $u\in Z^1\big(\pi_1(\Sigma_g),\mathfrak{h}\big)$,
\[u(ab)\ =\ u(a)+u(b),\ \forall\ a,b\in\mathbb{Z}\pi_1(\Sigma_g).\]
Hence, the values $\{u(a_i), u(b_i)\}_{i=1,...,g}$ determine $u$. Fixing $\rho_0$, any $\rho_1$ in $\mathcal{Y}_0$ is
\[\rho_1\ =\ \exp{\big(u\otimes \varepsilon+\mathcal{O}(\varepsilon^3)\big)}\cdot\rho_0,\]
for some $u\in Z^1\big(\pi_1(\Sigma_g),\mathfrak{h}\big)$. Therefore, an analytic chart of $\mathcal{Y}_0$ is $(B_{\varepsilon}(\rho_0), \log_{\varepsilon})$, where
\[\log_{\varepsilon}: \ B_{\varepsilon}(\rho_0)\ \rightarrow\ \mathbb{R}^{2n},\]
so that $\log_{\varepsilon}(\rho_0)=0$ and 
\[\log_{\varepsilon}(\rho)=\big(u(a_1),u(b_1),...,u(a_g),u(b_g)\big),\ \forall\ \rho_1\in B_{\varepsilon}(\rho_0).\]
Define a vector field $X$ on $\mathcal{Y}_0$ as a map from the set of differentiable functions to the set of functions on $\mathcal{Y}_0$. Let $f$ be a differentiable function on $\mathrm{U}(1)^{2g}$, and $X(\rho_0)\in Z^1\big(\pi_1(\Sigma_g),\mathfrak{h}\big)$ so that
\[X(\rho_0)\ =\ \big(u(a_1),u(b_1),\cdots, u(a_g),u(b_g)\big).\]
Then, 
\begin{equation}
    \label{DerivativeOfXOnf}
    \begin{aligned}
\big(X(f)\big)_{\rho_0}\ &=\ \frac{d}{dt}|_{t=0}f\big(\exp{\big(tX(\rho_0)\big)}\cdot \rho_0\big)\\
&=\ \frac{d}{dt}|_{t=0}f\big(\exp{tu(a_1)}\cdot\rho_0(a_1), \exp{tu(b_1)}\cdot \rho_0(b_1),\cdots, \exp{tu(a_g)}\cdot\rho_0(a_g), \exp{tu(b_g)}\cdot\rho_0(b_g)\big)\\
&=\ \sum_{i=1}^g\frac{\partial f}{\partial x_i}(\rho)\cdot \big(R_{\rho(a_i)}u(a_i)+R_{\rho(b_i)}u(b_i)\big),
\end{aligned}
\end{equation}
where $\frac{\partial f}{\partial x_i}$ is the derivative of $f$ relative to $i$-th component.

\begin{definition}
   \textit{ A vector field $X$ of $\mathcal{Y}_0$ is an assignment of $\rho_0$ to a 1-cocycle $Z^1\big(\pi_1(\Sigma_g),\mathfrak{h}\big)$ that is left invariant 
   \[X:\ \mathcal{Y}_0\ \rightarrow\ \mathfrak{h}^{2g}.\]   
   }
\end{definition}

\begin{proposition}
    If the vector field $X$ satisfies 
\[\mathcal{D}_{X(\rho_1)}X(\rho)\ =\ 0,\]
for any successive deformation $\rho$ and $\rho_s$, then $X$ is analytic. In particular, on $\mathcal{Y}_0$, an open ball $B_{\varepsilon}(\rho_0)$ centered at $\rho_0$ with radius small enough $\varepsilon$ consists of representations 
\[\rho_1\ =\ \exp{\big(u\otimes \varepsilon+\mathcal{O}(\varepsilon^3)\big)}\cdot \rho_0\]
so that $\rho_1$ is on $\mathcal{Y}_0$. A local flow of $X$ is an analytic mapping 
    \[\begin{aligned}
        \varphi:\ B_{\varepsilon}(\rho_0)\times(-\delta,\delta)\ &\rightarrow\ \mathcal{Y}_0,\\
        (\rho_1, t)\ &\mapsto\ \exp{\big(X(\rho_1)\otimes t+\mathcal{O}(t^2)\big)}\cdot\rho_0
    \end{aligned}\]
    for some $\delta>0$, so that for any $\rho_1\in B_{\varepsilon}(\rho_0)$, $\frac{\partial\varphi}{\partial t}=X\big(\varphi(\rho_1,t)\big)$.
\end{proposition}

\subsection{A vector bundle over the $\mathcal{Y}_0$}
Identify $\mathcal{Y}$ with $\mathrm{U}(1)^{2g}$, then $\mathcal{Y}_0$ is open in $\mathrm{U}(1)^{2g}$. Consider the natural projection,
\[\pi_N:\ \mathcal{N}\ \rightarrow\ \mathcal{Y}_0,\]
so that 
\[\pi_N^{-1}(\rho_0)\ =\ \mathcal{N}_{\rho_0}\ =\ Z^1\big(\pi_1(\Sigma_g),\mathfrak{h}_{\mathsf{Ad}_{\rho_0}}^{\perp}\big).\]
Identify the Lie algebra $\mathfrak{g}$ with the right-invariant vector fields of its Lie group $G$.
The space of 1-cocycles $ Z^1\big(\pi_1(\Sigma_g),\mathfrak{h}_{\mathsf{Ad}_{\rho_0}}^{\perp}\big)$ is a subspace of $(\mathfrak{h}^{\perp})^{2g}$. Use the Fox derivative (\cite{BillSymplecticNature}), a 1-cocycle $u\ =\ \big(u(a_1),u(b_1),...,u(a_g),u(b_g)\big)$ satisfies
\[\sum_{j=1}^{g}\big(\mathrm{Ad}_{\rho_0(\partial_{a_j}R)}u(a_j)+\mathrm{Ad}_{\rho_0(\partial_{b_j}R)}u(b_j)\big)\ =\ 0.\]
Let $\rho_0\in\mathcal{Y}_0$, and $B_{\varepsilon} (\rho_0)$ is the open ball centered at $\rho_0$ in $\mathcal{Y}_0$. Let $\rho\in B_{\varepsilon} (\rho_0)$, and $u$ be the 1-cocycle in  $Z^1\big(\pi_1(\Sigma_g),\mathfrak{h}\big)$ so that
\[\rho\ =\ \exp{\big(\varepsilon u+\mathcal{O}(\varepsilon^2)\big)}\cdot \rho_0.\]
Define a map
\begin{equation}
    \label{LocalTrivializationOfNormalBundle}
    \begin{aligned}
    \psi_{\rho_0}:\  \pi_N^{-1}\big(B_{\varepsilon}(\rho_0)\big)\ &\rightarrow\ B_{\varepsilon}(\rho_0) \times Z^1\big(\pi_1(\Sigma_g),\mathfrak{h}_{\mathsf{Ad}_{\rho_0}}^{\perp}\big)\\
    \tau\ &\mapsto\ \big(p_1(\tau),u(\tau)\big),
\end{aligned}
\end{equation}
where $p_1$ is the projection on the first component. Denote $p_1(\tau)=\rho$. Then
\[R_{\rho_0}\exp{(\varepsilon u)}\ =\ R_{\rho}\exp{0},\]
where $u\in T_{\rho_0}\mathcal{Y}_0$ and $0$ is the zero vector of $\mathfrak{h}$. Let 
\[p_2: \pi_N^{-1}\big(B_{\varepsilon}(\rho_0)\big)\ \rightarrow\ Z^1\big(\pi_1(\Sigma_g),\mathfrak{h}^{\perp}_{\mathrm{Ad}_{p_1(\tau)}}\big).\]
The 1-cocycle $u(\tau)$ satisfies the equation:
\[R_{\rho}u_1\ =\ R_{\rho_0}u(\tau),\]
where $u_1=p_2(\tau)$. Equivalently,
\[R_{\rho_0}R_{\exp{(\varepsilon u)}}p_2(\tau)\ =\ R_{\rho_0}R_{\exp{(\varepsilon u)}}u_1\ =\ R_{\rho_0}u(\tau),\]
where $R_{\exp{(\varepsilon u)}}=(\mathrm{Id}-\frac{\varepsilon}{2}\mathrm{ad}_{u})$. Hence,
\begin{equation}
    \label{LocalTrivializationOfNormalBundleN}
    \begin{aligned}
   (R_{\exp{(\varepsilon u)}}\circ p_2)(\tau) 
   =\  u_1-\frac{\varepsilon}{2}[u,u_1]+\mathcal{O}(\varepsilon^2)
   =\ u(\tau).
\end{aligned}
\end{equation}
In particular, 
\[\psi_{\rho_0}\ = \ \big(p_1,R_{\rho_0^{-1}\circ p_1}\circ p_2\big).\]
is analytic over $B_{\varepsilon}(\rho_0)$. Now, suppose $\rho$ is in $B_{\varepsilon}(\rho_0)\cap B_{\delta}(\rho_0)$, so that
\[\rho\ =\ \exp\big(u_0\otimes\varepsilon+\mathcal{O}(\varepsilon^2)\big)\cdot \rho_0\ =\ \exp\big(u_1\otimes\delta+\mathcal{O}(\delta^2)\big)\cdot \rho_1.\]
Therefore, the transition function is constant:
\[R_{\rho_1\rho_0^{-1}}:\ B_{\varepsilon}(\rho_0)\cap B_{\delta}(\rho_0)\ \rightarrow\ \mathrm{U}(1)^{2g}.\]

\subsection{A flat connection of the normal vector bundle}
Let $\sigma$ be a section of the vector bundle $\mathcal{N}$,
\[\sigma:\ B_{\varepsilon}(\rho_0)\ \rightarrow\ \pi_N^{-1}\big( B_{\varepsilon}(\rho_0)\big),\]
and $\sigma(\rho)=(\rho, u_0)$, for some $\rho\in B_{\varepsilon}(\rho_0)$, and 
\[\rho\ =\ \exp{\big(u\otimes \varepsilon+\mathcal{O}(\varepsilon^2)\big)}\cdot\rho_0.\]
Assume $u_1\in Z^1\big(\pi_1(\Sigma_g),\mathfrak{h}_{\mathrm{Ad}_{\rho}}^{\perp}\big)$ so that
\[R_{\rho}u_1\ =\ R_{\rho_0}u_0.\]
Hence, $u_0=R_{\exp{(\varepsilon u)}}u_1$. Now, let $X$ be a vector field of $\mathcal{Y}_0$. The multiplication of the group elements
\[\exp{\big(u_0\otimes r+v_0\otimes r^2+\mathcal{O}(r^3)\big)}\cdot\exp{\big(X(\rho_0)\otimes t+\mathcal{O}(t^3)\big)}\]
determines an infinitesimal deformation of $\rho_0$. Apply the BCHD formula to the multiplication of the exponential maps, 
\[\begin{aligned}
    u'(r,t)\ :=\ &u_1r+X(\rho_0)t+\frac{1}{2}[u_1, X(\rho_0)]rt+v_1 r^2+X(\rho_0)t^2+\frac{1}{2}[v_1,X(\rho_0)]r^2t+\mathcal{O}(t^3)+\mathcal{O}(r^3).
\end{aligned}\]
Let $X(\rho_0)=u$ in $Z^1\big(\pi_1(\Sigma_g),\mathfrak{h}\big)$. The bracket $[u_1,u]$ is a 2-cocycle in $Z^2\big(\pi_1(\Sigma_g),\mathfrak{h}^{\perp}_{\mathrm{Ad}_{}\rho_0}\big)$:
\[[u_1,u]\ =\ [u_1\wedge u]_{\mathfrak{su}(2)}.\]
Take the derivative relative to the parameter $s$ at $s=0$,
\begin{equation}
    \label{InfinitesimalSuccessiveDeformation}
     u_1+\frac{t}{2}[u_1, u]-\frac{t^2}{12}[u, [u, u_1]].
\end{equation}
Suppose $u$ translates $u_1$ by the infinitesimal transformation $d\varphi_{u}$ at $\rho_0$, where $d\varphi_{u}$ is an linear map(\ref{InfinitesimalSuccessiveDeformation}). Hence,
\[\begin{aligned}
    (d\varphi_{u})u_1\ &=\ u_1+\frac{t}{2}[u_1, u]+\frac{t^2}{2}[u_1,v]-\frac{t^2}{12}[u, [u, u_1]].
\end{aligned}\]
Hence, the variation of $u_1$ in the direction of $X(\rho_0)$ is
\[(\mathcal{D}_{u}u_1)_{\rho_0}\ :=\ \lim_{t\rightarrow0}\frac{1}{t}\big(u_1-(d\varphi_{u})u_1\big)\ =\ \frac{1}{2}[u_1,u].\]
Infinitesimally, the variation of $u$ along $X(\rho_0)$ is
\[\frac{d^2}{dtds}|_{t, s=0}u'(s,t)\ =\ \frac{1}{2}[X(\rho_0), u_1].\]
If $\mathcal{D}_{X(\rho_0)}u_1=0$, then the successive deformation is a deformation of $\rho_0$. Define a operator 
\[\mathcal{D}:\ \mathfrak{X}(\mathcal{Y}_0)\times\mathcal{N}_{\rho_0}\ \rightarrow\ \mathcal{N}_{\rho_0}\oplus\wedge^2\mathcal{N}_{\rho_0},\]
so that $\mathcal{D}_{X(\rho_0)}s:=[X(\rho_0),u_1]$. Hence, $\mathcal{D}$ is functional linear on the first component, and has degree $1$. \\\\
From the formula of vector field $X$ (\ref{DerivativeOfXOnf}), the differential of a function $f$ on $\mathcal{Y}_0$ is 
\[df\ =\ \sum_{i=1}^g\big(\frac{\partial f}{\partial x_i}R_{x_i}dx_i+\frac{\partial f}{\partial y_i}R_{y_i}dy_i\big).\]
Let 
\begin{itemize}
    \item $\mathcal{A}^0\big(\mathcal{Y}_0,\mathcal{N}\big)$ be the space of analytic functions over $\mathcal{Y}_0$ with values in $\mathcal{N}$;
    \item $\mathcal{A}^1\big(\mathcal{Y}_0, \mathcal{N}\oplus \wedge^2\mathcal{N}\big)$ be the space of 1-forms of $\mathcal{Y}_0$ with values in the bundle $\mathcal{N}\oplus \wedge^2\mathcal{N}$;
    \item $\mathcal{A}^2\big(\mathcal{Y}_0, \mathcal{N}\oplus \wedge^2\mathcal{N}\oplus\wedge^3\mathcal{N}\big)$ be the space of 2-forms of $\mathcal{Y}_0$ with values in the bundle $\mathcal{N}\oplus \wedge^2\mathcal{N}\oplus\wedge^3\mathcal{N}$;
\end{itemize}
Define 
\[\begin{aligned}
    \nabla:\ \mathcal{A}^0\big(\mathcal{Y}_0,\mathcal{N}\big)\ &\rightarrow\ \mathcal{A}^1\big(\mathcal{Y}_0, \mathcal{N}\oplus \wedge^2\mathcal{N}\big)\\
    fs\ &\mapsto\ \nabla(fs), 
\end{aligned}\] 
so that
\[\begin{aligned}
    \nabla(fs)\ &=\ (df)\otimes (s,0)+f\cdot (0,\mathcal{D}s),
\end{aligned}\]
where evaluating $\mathcal{D}$ on $X(\rho_0)$ is $\mathcal{D}_{X(\rho_0)}s$. 

\begin{theorem}
   \textit{The operator $\nabla=d+\mathcal{D}$ is a connection of the vector bundle $\mathcal{N}\rightarrow\mathcal{Y}_0$. In particular, the curvature $\mathcal{R}$ relative to $\nabla$ is $\mathcal{R}\ =\ \nabla\circ\nabla$, and $\mathcal{R}\equiv0$.}
\end{theorem}
\begin{proof}
$\nabla$ satisfies the following properties:
\begin{enumerate}
    \item 
    \[\begin{aligned}
        \nabla_{fX+gY}s\ &=\  \mathcal{D}_{fX+gY}s\\ 
        &=\ \frac{1}{2}[fX+gY, s]\\
        &=\ \frac{1}{2}f[X, s]+\frac{1}{2}g[Y,s]\\
        &=\ f\mathcal{D}_Xs+g\mathcal{D}_Ys\\
        &=\ f\nabla_Xs+g\nabla_Ys.
    \end{aligned}\]
    \item 
    \[\begin{aligned}
        \nabla_X(s_1+s_2)\ &=\ [X, s_1+s_2]\\
        &=\ \mathcal{D}_Xs_1+\mathcal{D}_Xs_2\\
        &=\ \nabla_Xs_1+\nabla_Xs_2.
    \end{aligned}\]
    \item The composition of $\nabla\circ\nabla$ is
    \[\mathcal{A}^0\big(\mathcal{Y}_0,\mathcal{N}\big)\ \xrightarrow{\nabla}\ \mathcal{A}^1(\mathcal{Y}_0, \mathcal{N}\oplus \wedge^2\mathcal{N}\big)\ \xrightarrow{\nabla}\ \mathcal{A}^2\big(\mathcal{Y}_0, \mathcal{N}\oplus \wedge^2\mathcal{N}\oplus\wedge^3\mathcal{N}\big)\]

    For any $f\in\mathcal{A}^0(\mathcal{Y}_0)$, and $s$ a section of the vector bundle $\mathcal{N}\rightarrow\mathcal{Y}_0$,
   \[\begin{aligned}
       \mathcal{R}(fs)\ &=\ \nabla\big((df)\otimes (s,0)+f\cdot (0,\mathcal{D}s)\big)\\ 
       &=\ (-1)^1\big(df\otimes(0,\mathcal{D}s)\big)+df\otimes (0,\mathcal{D}s)+ f\cdot(0, \mathcal{D}^2s)\\ 
       &=\ 0
   \end{aligned}\]
\end{enumerate}
\end{proof}

\subsection{A presymplectic embedding}
The fiber of the vector bundle $\mathcal{N}$ is 
\[\mathcal{N}_{\rho_0}\ =\ Z^1\big(\Gamma_g,\mathfrak{h}^{\perp}_{\mathrm{Ad}_{\rho_0}}\big),\]
which is presymplectic (\cite{BillSymplecticNature}). Let $\rho_0$ be a reducible noncentral representation, and $\mathrm{ad}E_{\rho_0}$ be the flat assoicated vector bundle $P_{\rho_0}\times_{\mathrm{U}(1)}\mathfrak{h}^{\perp}$. In terms of the deRham cohomolocy, let $\alpha_1$ and $\alpha_2$ be closed $1$-forms in $\Omega^1\big(\Sigma_g,\mathrm{ad}E_{\rho_0}\big)$, the presymplectic form in the fiber $\mathcal{N}_{\rho_0}$ is
    \[\omega^N_{\rho_0}(\alpha_1,\alpha_2)\ =\ \int_{\Sigma_g}\mathcal{B}_{\mathfrak{h}^{\perp}}(\alpha_1\wedge\alpha_2).\]
Still denote the presymplectic form of the vector space $Z^1\big(\Gamma_g,\mathfrak{h}^{\perp}_{\mathrm{Ad}_{\rho_0}}\big)$ in terms of the group cohomology as $\omega^N_{\rho_0}$, which has the explicit formula with the Fox derivative (\cite{BillSymplecticNature}).  Hence, the stratum $\big(\mathcal{Y}_0,\omega_{\rho_0}^{\mathcal{Y}_0}\big)$ is symplectic, where $\omega_{\rho_0}^{\mathcal{Y}_0}$ is the canonical Narasimhan-Atiyah-Bott-Goldman presymplectic form of $\mathsf{Hom}\big(\pi_1(\Sigma_g),\mathrm{U}(1)\big)$.

\begin{theorem}
    The vector bundle $\mathcal{N}\rightarrow\mathcal{Y}_0$ is a presymplectic vector bundle $(\mathcal{N},\omega^N)$.
\end{theorem}
\begin{proof}
    The fiber of the vector bundle $\mathcal{N}$ is presymplectic, and the presymplectic form $\omega^N$ is parallel along $\mathcal{Y}_0$. Let $s_1$ and $s_2$ be two sections of the bundle $\mathcal{N}$,
    \[\begin{aligned}
        \mathcal{D}_{X(\rho_0)}\omega^N(s_1, s_2)\ &=\ \mathcal{D}_{X(\rho_0)}\big(\omega^N(s_1,s_2)\big)+\omega^N\big(\mathcal{D}_{X(\rho_0)}s_1, s_2\big)+\omega^N\big(s_1, \mathcal{D}_{X(\rho_0)}s_2\big)\\
        &=\ \mathcal{D}_{X(\rho_0)}\big(\omega^N(s_1,s_2)\big).\\
     &=\ \omega^N\big(d\varphi_{X(\rho_0)}s_1, d\varphi_{X(\rho_0)}s_2\big)\\
     &=\ \omega^N\big([X(\rho_0),s_1], [X(\rho_0),s_2]\big)\\
    &=\ 0.
    \end{aligned} \]
    Finally, $\omega^N$ is closed under $\nabla$.\\\\
    Since the kernel of $\omega^{\mathcal{Y}_0}$ is trivial, the constant-rank embedding theorem (\cite{SjamaarStratified}) implies that the Whitney sum determines the presymplectic structure on the reducible noncentral stratum.

\end{proof}

\begin{theorem}
 \label{PresymplecticEmbedding}
\textit{
The presymplectic vector bundle $\big(\mathcal{N}, \omega^{N}\big)$ over the symplectic base $\big(\mathcal{Y}_0,\omega^{\mathcal{Y}_0}\big)$ determines a pre-symplectic embedding 
\[i:\ \big(\mathcal{Y}_0,\omega^{\mathcal{Y}_0}\big)\ \hookrightarrow\ \Big(\mathsf{Hom}\big(\pi_1(\Sigma_g),\mathrm{SU}(2)\big),\omega_{\mathsf{NABG}}\Big)\]
so that $i^*\omega_{\mathsf{NABG}}\ =\ \omega^{\mathcal{Y}_0}$. Therefore, the canonical Narasimhan-Atiyah-Bott-Goldman presymplectic form is 
\begin{equation}
\label{PresymplecticFormOnY0}
    \omega_{\mathsf{NABG}}\ =\ \omega^{\mathcal{Y}_0}+\omega^{N}.
\end{equation}}
\end{theorem}

\subsection{The cone bundle over $\mathcal{Y}_0$}
On the stratum $\mathsf{Hom}\big(\pi_1(\Sigma_g),\mathrm{SU}(2)\big)_{\mathrm{U}(1)}$, at $\rho_0$, the tangent space is 
\[\mathcal{T}_{\rho_0}\ \cong\ Z^1\big(\pi_1(\Sigma_g),\mathfrak{su}(2)_{\mathrm{Ad}_{\rho_0}}\big)\ \cong\  Z^1\big(\pi_1(\Sigma_g),\mathfrak{h}\big)\times  Z^1\big(\pi_1(\Sigma_g),\mathfrak{h}^{\perp}_{\mathrm{Ad}_{\rho_0}}\big),\]
and 
\[Z_0\big(\pi_1(\Sigma_g),\mathfrak{h}_{\mathrm{Ad}_{\rho_0}}^{\perp}\big)\ = \ 0.\]
Define a pairing on $\mathcal{T}_{\rho_0}$:
\begin{equation}
\label{PairingOnTangentSpaceSplitingTangentialAndTransverseDirections}
    g_{\rho_0}:\ \mathcal{T}_{\rho_0}\times \mathcal{T}_{\rho_0}\ \xrightarrow{[\cdot\ \wedge\ \cdot]_{\mathfrak{su}(2)}}\ Z^2(\pi_1(\Sigma_g),\mathfrak{su}(2)_{\mathrm{Ad}_{\rho_0}}\big)\ \xrightarrow{\cap[\Sigma_g]}\ Z_0\big(\pi_1(\Sigma_g),\mathfrak{h}\big),
\end{equation}
where $[\cdot\ \wedge\ \cdot]_{\mathfrak{su}(2)}$ is the quadratic bracket on $1$-cocycles, and $[\Sigma_g]$ is the fundamental 2-cocycles of the closed surface $\Sigma_g$ in $H_2\big(\pi_1(\Sigma_g),\mathbb{R} \big)$. Using group cohomology and the Fox derivative, the fundamental 2-cocycle has a nice expression (\cite{BillSymplecticNature}). Indeed,
\[\begin{aligned}
    Z^1\big(\pi_1(\Sigma_g),\mathfrak{h}\big)\ &\cong\ \mathbb{R}^{2g}\\
    Z_0\big(\pi_1(\Sigma_g),\mathfrak{h}\big)\ &\cong\ \mathbb{R}.
\end{aligned}\]

\begin{theorem}
\label{DecompositionOFConeIntoConeInTheNormalBundle}
    \textit{Let $\rho_0$ be in $\mathcal{Y}_0$, the symmetric pairing $g_{\rho_0}$ determines a decomposition of the tangent quadratic cone $\mathcal{L}_{\rho_0}$ into the direct sum of a Euclidean space $\mathbb{R}^{2g}$ with a quadratic cone,
    \[\mathcal{L}_{\rho_0}\ \cong\ Z^1\big(\pi_1(\Sigma_g),\mathfrak{h}\big)\oplus \mathcal{C}_{\rho_0},\]
    where 
   \begin{equation}
       \label{QuadraicConeInNormalBundle}
       \mathcal{C}_{\rho_0}\ =\ \{u\in Z^1\big(\pi_1(\Sigma_g),\mathfrak{h}^{\perp}_{\mathrm{Ad}_{\rho_0}}\big):\ u\in\mathcal{L}_{\rho_0}\}.
   \end{equation}
    In particular, in the vector bundle $\mathcal{N}$, the quadratic cone $\mathcal{C}_{\rho_0}$ is a cone in the fiber $\mathcal{N}_{\rho_0}$.}
\end{theorem}

\begin{proof}
The pairing $g_{\rho_0}$ is symmetric, and decomposes the Zariski tangent space at $\rho_0$ to be
\[\mathcal{T}_{\rho_0}\ =\ Z^1\big(\pi_1(\Sigma_g),\mathfrak{h}\big)\oplus Z^1\big(\pi_1(\Sigma_g),\mathfrak{h}^{\perp}_{\mathrm{Ad}_{\rho_0}}\big).\]
Let $u_1$ and $u_2$ be in $Z^1\big(\pi_1(\Sigma_g),\mathfrak{su}(2)_{\mathrm{Ad}_{\rho_0}}\big)$. 
\begin{itemize}
    \item If one of the 1-cocycles is in $Z^1\big(\pi_1(\Sigma_g),\mathfrak{h}\big)$, then $g_{\rho_0}(u_1,u_2)$ vanishes; 
    \item If one of the two components is a 1-coboundary, suppose $u_1$ is in $B^1\big(\pi_1(\Sigma_g),\mathfrak{h}_{\mathrm{Ad}_{\rho_0}}^{\perp}\big)$. Then there exists a $u_2$ in $Z^1\big(\pi_1(\Sigma_g),\mathfrak{h}^{\perp}_{\mathrm{Ad}_{\rho_0}}\big)$, so that $[u_1,u_2]\in Z^2\big(\pi_1(\Sigma_g),\mathfrak{h}\big)$ and $g_{\rho_0}(u_1,u_2)\neq0$. 
\end{itemize}
Therefore,
\[Z^1\big(\pi_1(\Sigma_g),\mathfrak{h}\big)=\ker{g}_{\rho_0},\]
and $g_{\rho_0}$ is nondegerate on $Z^1\big(\pi_1(\Sigma_g),\mathfrak{h}^{\perp}_{\mathrm{Ad}_{\rho_0}}\big)$.
\end{proof}

\begin{theorem}
\label{QuadraticConeFiberBundle}
    \textit{The restriction of $g$ on $\mathcal{N}$ is symmetric nondegenerate. Therefore, the fiber bundle 
    \[\pi_C:\ \mathcal{C}\rightarrow\mathcal{Y}_0\]
    is a subbundle of $\mathcal{N}\rightarrow\mathcal{Y}_0$ so that $\pi_C^{-1}(\rho_0)=\mathcal{C}_{\rho_0}$,} 
    \[\begin{tikzcd}
	{\mathcal{C}} & {\mathcal{N}} \\
	& {\mathcal{Y}_0}
	\arrow[from=1-1, to=1-2]
	\arrow[from=1-1, to=2-2]
	\arrow[from=1-2, to=2-2]
\end{tikzcd}\]
\end{theorem}
\begin{proof}
 The vector bundle $\mathcal{N}$ admits a fiberwise nondegenerate symmetric $g_{\rho_0}$.\\\\ 
 Furthermore, the bilinear two-form $g$ is parallel relative to the connection $\nabla$. Let $X$ be a vector field over $\mathcal{Y}_0$, then 
\[\nabla_{X(\rho_0)}g(s_1,s_2)\ =\ g\big(\nabla_{X(\rho_0)}s_1,s_2\big)+g\big(s_1,\nabla_{X(\rho_0)}s_2\big).\]
Since the graded Jacobian is
\[\begin{aligned}
   0\ &=\ [[X,s_1],s_2]
+[[s_1,s_2],X]+[[s_2,X],s_1]\\
&=\ [[X,s_1],s_2]
-[X,[s_1,s_2]]+[[X,s_2],s_1].
\end{aligned}\]
Eventually, the local trivialization of the bundle $\mathcal{N}$ preserves the tangent quadratic cone:
\[[R_{\exp{\varepsilon u}}u_1,R_{\exp{\varepsilon u}}]\ =\ 0,\]
where $R_{\exp{\varepsilon u}}$ determines $\psi_{\rho_0}$ in the local trivialization of the bundle $\mathcal{N}$ at $\rho_0$ (\ref{LocalTrivializationOfNormalBundle}).

\end{proof}

\begin{theorem}
\label{TubularNeighborhoodTheorem}
    \textit{Let $\mathcal{U}_{\mathcal{R}}$ be the set of all infinitesimal deformations of representations in $\mathcal{Y}_0$,
    \[\{\exp{\big(u\otimes t+v\otimes t^2+\mathcal{O}(t^3)\big)}\cdot\rho_0:\ \forall\ \rho_0\in \mathcal{Y}_0,\ \text{and}\ \forall\  u\in\mathcal{C}_{\rho_0}\}.\]
    The analytic map 
    \[\exp:\ \mathcal{U}_0\ \rightarrow\ \mathcal{U}_{\mathcal{R}}\]
    determines an tubular neighborhood of $\mathcal{Y}_0$ in $\mathsf{Hom}\big(\pi_1(\Sigma_g),\mathrm{SU}(2)\big)$, where $\mathcal{U}_0$ is a neighborhood of $\mathcal{Y}_0$ in the bundle $\mathcal{N}$.}
\end{theorem}

\newpage
\section{The Symplectic Nature Of The Stratum Of Reducible Noncentral Representations in the $\mathrm{SU}(2)$ Character Variety}
The group $\mathrm{SU}(2)$ acts by adjoint action on the representation variety $\mathrm{Hom}\big(\pi_1(\Sigma_g),\mathrm{SU}(2)\big)$. The group $\mathrm{SU}(2)$ acts on $\mathrm{Hom}\big(\pi_1(\Sigma_g),\mathrm{SU}(2)\big)$ is not locally free at reducible noncentral representations. The isotropy group $\mathrm{U}(1)$ acts on the Zariski tangent space at a reducible noncentral representation nontrivially.\\\\  
Take the quotient of the adjoint $\mathrm{SU}(2)$ action, the $\mathrm{SU}(2)$-character variety is 
\begin{equation}
    \label{CharacterVarietyDefinition}
    \chi\big(\pi_1(\Sigma_g),\mathrm{SU}(2)\big)\ :=\ \mathrm{Hom}\big(\pi_1(\Sigma_g),\mathrm{SU}(2)\big)/ \mathrm{SU}(2).
\end{equation}
The reducible noncentral stratum in the character variety is 
\[\chi_0\ :=\ \widetilde{\mathcal{Y}}_0/\mathrm{SU}(2),\]
where $\widetilde{\mathcal{Y}}_0=\mathrm{SU}(2)\cdot\mathcal{Y}_0$. Indeed, the neighborhood $\mathcal{U}_0$ of the zero section in the normal vector bundle $\mathcal{N}$ provides a picture of the local structure of $\chi_0$, where 
\[T_{[\rho_0]}\chi_0\ \cong H^1\big(\pi_1(\Sigma_g),\mathfrak{su}(2)_{\mathrm{Ad}_{\rho_0}}\big)\big/\mathrm{U}(1),\]
is the quotient by the isotropy group action. In particular, the germ of $\chi_0$ at the class $[\rho_0]$ is isomorphic to $\widetilde{\mathcal{L}}_{\rho_0}/\mathrm{U}(1)$, the quotient of the reduced cone  $\widetilde{\mathcal{L}}_{\rho_0}$ which is in $H^1\big(\pi_1(\Sigma_g),\mathfrak{su}(2)_{\mathrm{Ad}_{\rho_0}}\big)$. The Hamiltonian space and symplectic reduction are the main tools for analyzing the symplectic structure of $\chi_0$. \\\\
The q-Hamiltonian space with group-valued moment map reveals the relationship between the symplectic structures of $\mathcal{Y}_0$ and $\chi_0$. Let $T\cong \mathrm{U}(1)$ be a maximal torus of $\mathrm{SU}(2)$. Its associated Weyl group is
\begin{equation}
    \label{WeylGroupDefinition}
    W\ :=\ N_{\mathrm{SU}(2)}(T)/T\ \cong\ \mathbb Z/2\mathbb Z,
\end{equation}
where $N_{\mathrm{SU}(2)}(T)$ represents the normalizer of $T$ in $\mathrm{SU}(2)$. Write $W=\{1,w\}$, where $w$ is the unique nontrivial element and satisfies $w^2=1$. The map
\begin{equation}
    \label{GroupValuedMomentMap}
    \begin{aligned}
    \mu:\ \mathrm{Hom}\big(\pi_1(\Sigma_g),\mathrm{SU}(2)\big)\ &\rightarrow\ \mathrm{SU}(2)\\
\rho\ &\mapsto\  \rho(a_1)\rho(b_1)\rho(\overline{a_1})\rho(\overline{b_1})\cdots\rho(a_g)\rho(b_g)\rho(\overline{a_g})\rho(\overline{b_g})
\end{aligned}
\end{equation}
is a constant map. In particular, the space $\big(\mathcal{Y}_0, \omega, \mathrm{Ad}_{W},\mu|_{\mathcal{Y}}\big)$ is quasi-Hamiltonian.

\subsection{The isotropy group action}
 Let $\rho_0$ be a reducible noncentral representation. While $\mathrm{SU}(2)$ acting on the representation variety $\mathsf{Hom}\big(\pi_1(\Sigma_g),\mathrm{SU}(2)\big)$ by conjugation, the subgroup $\mathrm{U}(1)$ centralizes $\rho_0$ but acts on the tangent space by conjugation. In particular, the action $\mathrm{Ad}_{\mathrm{U}(1)}$ preserves the quadratic cone $\mathcal{L}_{\rho_0}$.\\\\
  Let $h\in\mathrm{U}(1)$, 
\[\begin{aligned}
    \frac{d}{dt}\big|_{t=0}\mathrm{Ad}_h\rho_t\ &=\  \frac{d}{dt}\big|_{t=0}h\cdot \exp{\big(u\otimes t+ v\otimes t^2+\mathcal{O}(t^3)\big)}\cdot\rho_0\cdot h^{-1}\\
    &=\ \frac{d}{dt}\big|_{t=0}\exp{\big(\mathrm{Ad}_hu\otimes t+\mathrm{Ad}_{h}v\otimes t^2+\mathcal{O}(t^3)\big)}\cdot \rho_0.
\end{aligned}\]
Hence, the group $\mathrm{U}(1)$ acts on the Zariski tangent space $\mathcal{T}_{\rho_0}$ by $\mathrm{Ad}_{\mathrm{U}(1)}$. \\\\
Now, let $u\in \mathcal{L}_{\rho_0}$, and $v$ be the 1-cochain so that for all $x,y\in \Gamma_g^1$,
\[[u,u](x,y)\ =\ -2(\delta v)(x,y).\]
Applying $\mathrm{Ad}_h$ to the left, 
\[\begin{aligned}
    \mathrm{Ad}_h[u,u](x,y)\ &=\ \mathrm{Ad}_h[u(x),\mathrm{Ad}_{\rho_0(x)}u(y)]\\ 
    &=\ [\mathrm{Ad}_hu(x),\mathrm{Ad}_{\rho_0(x)}\mathrm{Ad}_hu(y)]\\ 
    &=\ [\mathrm{Ad}_hu,\mathrm{Ad}_hu](x,y)
\end{aligned}\]
Applying $\mathrm{Ad}_h$ to the right,
\[\mathrm{Ad}_h(\delta v)\ =\ \delta(\mathrm{Ad}_hv).\]
Therefore, $\mathrm{Ad}_{h}$ is an automorphism on $\mathcal{L}_{\rho_0}$, for all $h\in\mathrm{U}(1)$,
\begin{equation}
    \label{EquivariantIsotropyGroupActionOnCone}
     [\mathrm{Ad}_hu,\mathrm{Ad}_hu](x,y)=\delta(\mathrm{Ad}_hv),\ \forall\ h\in\mathrm{U}(1).
\end{equation}
 
\subsection{The local structure of the reducible noncentral stratum in the character variety}
The group $\mathrm{SU}(2)$ acts on $\mathsf{Hom}\big(\pi_1(\Sigma_g),\mathrm{SU}(2)\big)$ by conjugation. Let $\mathcal{O}_{\rho_0}$ be the orbit $\mathrm{SU}(2)\cdot\rho_0$, then $\mathcal{O}_{\rho_0}$ is closed in $\mathsf{Hom}\big(\pi_1(\Sigma_g),\mathrm{SU}(2)\big)$. The stratum $\widetilde{\mathcal{Y}}_0$ is the disjoint union of orbits
\[\widetilde{\mathcal{Y}}_0\ =\ \bigsqcup_{\rho_0\in\mathcal{Y}}\mathcal{O}_{\rho_0},\]
that is closed in the tubular neighborhood $\mathcal{U}_{\mathcal{R}}$(\ref{TubularNeighborhoodTheorem}) of $\mathcal{Y}_0$ in $\mathsf{Hom}\big(\pi_1(\Sigma_g),\mathrm{SU}(2)\big)$. Up to the symmetric pairing $g_{\rho_0}$ (\ref{PairingOnTangentSpaceSplitingTangentialAndTransverseDirections}) at $\rho_0$, the Zariski tangent space at $\rho_0$ is
\[\mathcal{T}_{\rho_0}\ =\ Z^1\big(\pi_1(\Sigma_g),\mathfrak{h}\big)\oplus Z^1\big(\pi_1(\Sigma_g),\mathfrak{h}_{\mathrm{Ad}_{\rho_0}}^{\perp}\big).\]
The Zariski tangent space of the orbit $\mathcal{O}_{\rho_0}$ at $\rho_0$ is isomorphic to the space of 1-coboundaries:
\[T_{\rho_0}\mathcal{O}_{\rho_0}\ \cong\ B^1\big(\Gamma_g, \mathfrak{h}_{\mathrm{Ad}_{\rho_0}}^{\perp}\big).\]
In particular, the presymplectic form $\omega_{\mathsf{NABG}}=\omega^{\mathcal{Y}_0}+\omega^N$ has the kernel $B^1\big(\Gamma_g, \mathfrak{h}_{\mathrm{Ad}_{\rho_0}}^{\perp}\big)$at $\rho_0$. Hence, the orbit $\mathcal{O}_{\rho_0}$ is isotropic in the tubular neighborhood $\mathcal{U}_{\mathcal{R}}$,
 \[\frac{\big(T_{\rho_0}\mathcal{O}_{\rho_0}\big)^{\omega_{\rho_0}}}{T_{\rho_0}\mathcal{O}_{\rho_0}}\ \cong\ Z^1\big(\Gamma_g,\mathfrak{h}\big)\oplus H^1\big(\Gamma_g,\mathfrak{h}_{\mathrm{Ad}_{\rho_0}}^{\perp}\big),\] 
where $\omega_{\mathrm{NABG}}$ is non degenerate and reduced to the symplectic form 
   \[\Omega_{\rho_0}\ :=\ \omega_{\rho_0}^{\mathcal{Y}_0}+\Omega^{N}_{\rho_0}.\]

\subsection{A symplectic slice of the reducible noncentral stratum}
 Identify $H^2\big(\pi_1(\Sigma_g),\mathfrak{h}\big)$ with  $\mathfrak{h}^*$, define 
\[\begin{aligned}
    \Psi_{\rho_0}:\ H^1\big(\pi_1(\Sigma_g),\mathfrak{su}(2)_{\mathrm{Ad}_{\rho_0}}\big)\ &\rightarrow\ \mathfrak{h}^*\\
    u_0\ &\mapsto\ [u_0, u_0],
\end{aligned}\]
then $\Psi_{\rho_0}$ is a moment map (The same proof as in \ref{SymplecticNormalConeInCharacterVariety}). The space
\[\Big(H^1\big(\pi_1(\Sigma_g),\mathfrak{su}(2)_{\mathrm{Ad}_{\rho_0}}\big), \Omega_{\rho_0}, \mathrm{Ad}_{\mathrm{U}(1)}, \Psi_{\rho_0}\Big)\]
is a Hamiltonian space. From the symplectic reduction, the space
\[\begin{aligned}
    \widetilde{\mathcal{L}}_{\rho_0}\ :=&\  \Psi_{\rho_0}^{-1}(0)/\mathrm{U}(1)\\
    =&\ \Big(Z^1\big(\pi_1(\Sigma_g),\mathfrak{h}\big)\oplus\frac{\mathcal{C}_{\rho_0}}{B^1\big(\pi_1(\Sigma_g),\mathfrak{h}_{\mathrm{Ad}_{\rho_0}^{\perp}}\big)}\Big)\big/\mathrm{U}(1)\\
    =&\ Z^1\big(\pi_1(\Sigma_g),\mathfrak{h}\big)\oplus\widetilde{\mathcal{C}}_{[\rho_0]},
\end{aligned}\]
is symplectic with the symplectic form $\widetilde{\Omega}_{\rho_0}$. Suppose the reduction is 
\[\mathrm{Red}:\ \Psi_{\rho_0}^{-1}(0)\ \rightarrow\ \widetilde{\mathcal{L}}_{\rho_0},\]
then
\begin{equation}
    \label{ReducedSymplecticFormOfTheGermInCharacterVariety}
    \mathrm{Red}^*\widetilde{\Omega}_{\rho_0}\ =\ \Omega_{\rho_0}.
\end{equation}
The closed surface $\Sigma_g$ is the Eilenberg–Mac Lane space of type $K\big(\pi_1(\Sigma_g),1\big)$, 
\[H^i\big(\pi_1(\Sigma_g),\mathfrak{g}_{\mathrm{Ad}_{\rho_0}}\big)\ \cong\ H^i_{dR}\big(\Sigma_g,\mathrm{ad}P_{\rho_0}\big), \ i=0,1,2,\]
where $\mathrm{ad}P_{\rho_0}$ is the flat vector bundle over $\Sigma_g$ with the connection corresponding to $\rho_0$. Hence, the space
\[\Big(H^1\big(\Sigma_g,\mathrm{ad}P_{\rho_0}\big), \Omega_{\rho_0}, \mathrm{Ad}_{\mathrm{U}(1)},  (\Psi_{\rho_0})_{dR}\Big)\]
is Hamiltonian, where
\[\begin{aligned}
    (\Psi_{\rho_0})_{dR}:\ H^1\big(\Sigma_g,\mathrm{ad}P_{\rho_0}\big)\ &\rightarrow\ H^2\big(\Sigma_g,\mathrm{ad}P_{\rho_0}\big)\\
    \tau\ &\mapsto\ [\tau\wedge\tau].
\end{aligned}\]
\begin{proposition}
    The space $\big( (\Psi_{\rho_0})_{dR}^{-1}(0),\ \Omega_{\rho_0},\ \mathrm{Ad}_{U}(1),\ \mu_{\rho_0}\big)$ is a quasi-Hamiltonian space, where 
    \[\mu_{\rho_0}\ =\ \exp\circ\int_{\Sigma_g}\ (\Psi_{\rho_0})_{dR}:\ (\Psi_{\rho_0})_{dR}^{-1}(0)\ \rightarrow\ \mathrm{U}(1).\]
\end{proposition}
\begin{proof}
    Let $\theta$ denote the $\mathfrak h$-valued Maurer–Cartan form on $\mathrm{U}(1)$, where $\mathfrak h=\mathfrak u(1)$, then $\theta\ =\ h^{-1}dh$, for $h\in\mathrm{U}(1)$. Hence,
\[\mu^*(\theta,\xi)\ =\ (\exp\circ\ \int_{\Sigma_g}(\Psi_{\rho_0})_{dR})^*(\theta,\xi).\]
Now, 
\[\langle\mu_{\rho_0}^*(\theta,\xi),u\rangle\ =\ \langle\theta\big(d\mu_{\rho_0}(u)\big),\xi\rangle\ =\ \langle\int_{\Sigma_g}[u,u],\xi\rangle\ =\ \int_{\Sigma_g}\mathcal{B}([u,u],\xi),\]
where $\langle\ \cdot,\cdot\ \rangle$ is the Killing form of $\mathrm{su}(2)$. On the other hand,
\[\iota_{u_{\xi}}\Omega_{\rho_0}(u)\ =\ \Omega_{\rho_0}([\xi,u],u)\ =\ \int_{\Sigma_g}\mathcal{B}([\xi,u]\wedge u)\ =\ \frac{1}{2}\int_{\Sigma_g}\mathcal{B}([u\wedge u],\xi).\]
Hence, 
\[\iota_{u_{\xi}}\Omega_{\rho_0}\ =\ \frac{1}{2}\mu_{\rho_0}^*(\theta,\xi) \]
In particular, the 2-form $\Omega_{\rho_0}$ is nondegenerate on $(\Psi_{\rho_0})_{dR}^{-1}(0)$ and is closed. Therefore, the space $\big( (\Psi_{\rho_0})_{dR}^{-1}(0),\ \Omega_{\rho_0},\ \mathrm{Ad}_{U}(1),\ \mu_{\rho_0}\big)$ is q-Hamiltonian (\cite{Alekseev1997LieGV}).

\end{proof}

\begin{proposition} 
\label{SymplecticSlice}
A neighborhood of the orbit $\mathcal{O}_{\rho_0}$ is equivariantly homeomorphic to a neighborhood of the zero section of
\[\mathcal{Z}_{\rho_0}\ =\ \mathrm{SU}(2)\times_{\mathrm{U}(1)}\Big(\mathrm{SU}(2)/\mathrm{U}(1)\times(\Psi_{\rho_0})_{dR}^{-1}(0)\Big),\]
 In particular, the embedding $\mathrm{SU}(2)/\mathrm{U}(1)\rightarrow \mathcal{Z}_{\rho_0}$ is isotropic and the normal bundle is 
\[\mathrm{SU}(2){\times}_{\mathrm{U}(1)}(\Psi_{\rho_0})_{dR}^{-1}(0).\]
\end{proposition}
\begin{proof}
Let $D\big(\mathrm{SU}(2)\big)$ be the double Lie group 
\[\mathrm{SU}(2)\times \mathrm{SU}(2).\]
Indeed, $D\big(\mathrm{SU}(2)\big)$ is a q-Hamiltonian space under the action of $\mathrm{SU}(2)\times\mathrm{SU}(2)$,
\[(g_1,g_2)\cdot (a,b)\ :=\ (g_1ag_2^{-1}, g_2bg_1^{-1}).\]
The corresponding group-valued moment map is
\[(\mu_1,\mu_2)(a,b)\ =\ (ab, a^{-1}b^{-1}).\]
Now, under the Hopf fibration, 
\[\mathrm{U}(1)\ \rightarrow\ \mathrm{SU}(2)\ \xrightarrow{p}\ \mathrm{SU}(2)/\mathrm{U}(1),\]
the section $s$ is
\[\begin{aligned}
    s: \mathrm{SU}(2)/\mathrm{U}(1)\ &\rightarrow\ \mathrm{SU}(2)\\
    [x]\ &\mapsto\ s([x]).
\end{aligned}\]
in an open neighborhood of the identity coset, $e\mathrm{U}(1)$ so that $p(s([x]))=[x]$. Embed $\mathrm{U}(1)$ into $D\big(\mathrm{SU}(2)\big)$ by the composition:
\[\begin{aligned}
    \mathrm{U}(1)\times\{1\}\ &\hookrightarrow\ \mathrm{SU}(2)\times\mathrm{SU}(2)/\mathrm{U}(1)\ \subset\ D\big(\mathrm{SU}(2)\big)\\
    (h,1)\ &\mapsto\big(h,e\mathrm{U}(1)\big).
\end{aligned}\]
Suppose $u\in (\Psi_{\rho_0})_{dR}^{-1}(0)$, and $(h',1)$ acts on the product of $D\big(\mathrm{SU}(2)\big)$ and $(\Psi_{\rho_0})_{dR}^{-1}(0)$ simultaneously:
\[(h', 1)\cdot \big(h,e\mathrm{U}(1), u\big)\ =\ (h'h, e\mathrm{U}(1)\cdot h^{-1}, h\cdot u)\ =\big(h,e\mathrm{U}(1), u\big)\ =\ (h'h, e\mathrm{U}(1), \mathrm{Ad}_h u).\]
Reduce the $\mathrm{U}(1)$-action results in the iterated bundle,
\[\mathrm{SU}(2)\times \big(\mathrm{SU}(2)/\mathrm{U}(1)\times (\Psi_{\rho_0})_{dR}^{-1}(0)\big)\ \rightarrow\ \mathrm{SU}(2)\times_{\mathrm{U}(1)}\big(\mathrm{SU}(2)/\mathrm{U(1)}\times (\Psi_{\rho_0})_{dR}^{-1}(0)\big).\]
Hence, the group $\mathrm{U}(1)$ acts on the left in an equivariant way, hence
\[\mathcal{Z}_{\rho_0}\ :=\ \mathrm{SU}(2)\times_{\mathrm{U}(1)}\big(\mathrm{SU}(2)/\mathrm{U(1)}\times (\Psi_{\rho_0})_{dR}^{-1}(0)\big),\]
so that the zero section of $\mathcal{Z}_{\rho_0}$ is $\mathrm{SU}(2)/\mathrm{U}(1)$ which is homeomorphic to the orbit $\mathrm{SU}(2)\cdot\rho_0$.

\end{proof}

\begin{corollary}
    The symplectic slice at $[\rho_0]$ in the reducible noncentral stratum $\chi_0$ is 
    \[\Psi_{[\rho_0]}^{-1}(0)/\mathrm{U}(1)\ =\ \widetilde{\mathcal{L}}_{[\rho_0]},\]
    with the symplectic form $\widetilde{\Omega}_{[\rho_0]}$.
\end{corollary}

\subsection{A symplectic quadratic cone fiber bundle over the Kummer variety}
In the presymplectic vector bundle $\mathcal{N}\rightarrow\mathcal{Y}_0$,
\[\begin{aligned}
    T_{\rho_0}\mathcal{Y}_0\ &=\ Z^1\big(\pi_1(\Sigma_g),\mathfrak{h}\big)\\
    \mathcal{N}_{\rho_0}\ &=\ Z^1\big(\pi_1(\Sigma_g),\mathfrak{h}_{\mathrm{Ad}_{\rho_0}}^{\perp}\big).
\end{aligned}\]
As discussed above, the kernel of $\omega^{N}$ at representation $\rho_0$ is isomorphic to the tangent space of the orbit space
\[T_{\rho_0}\big(\mathrm{SU}(2)\cdot\rho_0\big)\ \cong\ B^1\big(\pi_1(\Sigma_g),\mathfrak{h}_{\mathrm{Ad}_{\rho_0}}^{\perp}\big).\]
Hence, the $\mathrm{SU}(2)$-action determines a bundle automorphism
\[\mathcal{N}\ \rightarrow\ \mathcal{N},\]
which covers the identity on the tangent direction. \\\\
Let $W$ be the Weyl group associated with the maximal torus of $\mathrm{U}(1)$ in $\mathrm{SU}(2)$, and $w$ is the nontrivial element in $W$ (\ref{WeylGroupDefinition}). Indeed, $w$ identifies $\rho_0$ with $\rho_0^{-1}$ on $\mathcal{Y}_0$. The Kummer variety $\mathcal{K}$ is
\[\mathcal{K}\ =\ \chi\big(\pi_1(\Sigma_g),\mathrm{U}(1)\big)/\{\pm1\}.\]
The $\mathrm{SU}(2)$-action on the stratum $\mathcal{Y}_0$ in 
$\mathsf{Hom}\big(\pi_1(\Sigma_g),\mathrm{SU}(2)\big)$ determines $\mathcal{K}_0$, the noncentral representation classes in $\mathcal{K}$.

\begin{theorem}
\label{Quotient Bundle Over The Kummer Variety}
The vector bundle $\mathcal{N}\rightarrow\mathcal{Y}_0$ is $\mathrm{SU}(2)$-equivariant, and the quotient $\Big(\mathcal{N}/\mathrm{SU}(2),\widetilde{\Omega}^{N}\Big)$ is a symplectic bundle over the noncentral classes in the Kummer variety $\mathcal{K}_0$. The fiber over $[\rho_0]$ is isomorphic to 
\[\Big(H^1\big(\pi_1(\Sigma_g),\mathfrak{h}_{\mathrm{Ad}_{\rho_0}}^{\perp}\big),\widetilde{\Omega}_{[\rho_0]}^N\Big).\]
\end{theorem}

\begin{proof}
Infinitesimally, the $\mathrm{SU}(2)$-action on $\mathcal{N}$ over the coordinate chart $\mathcal{B}_{\varepsilon}(\rho_0)$ is
\[\begin{aligned}
    \exp{t\xi}\cdot\psi_{\rho_0}:\  \pi_N^{-1}\big(B_{\varepsilon}(\rho_0)\big)\ &\rightarrow\ B_{\varepsilon}(\rho_0) \times Z^1\big(\pi_1(\Sigma_g),\mathfrak{h}_{\mathsf{Ad}_{\rho_0}}^{\perp}\big)\\
    \exp{t\xi}\cdot\tau\ &\mapsto\ \big(p_1(\tau),\xi(\tau)\big),
\end{aligned}\]
where $p_1(\tau)=\exp{\big(\varepsilon u +\mathcal{O}(\varepsilon^2)\big)}\cdot\rho_0$, and $\xi\in\mathfrak{su}(2)$. From the local trivialization of the bundle $\mathcal{N}$ (\ref{LocalTrivializationOfNormalBundleN}),
\[\begin{aligned}
    \xi(\tau)\ &=\ \varepsilon u+ u_1-\frac{\varepsilon}{2}\mathrm{ad}_{u}u_1+ (\xi-\mathrm{Ad}_{\rho_0}\xi)\\
    &=\ u(\tau)+(\xi-\mathrm{Ad}_{\rho_0}\xi)
\end{aligned}.\]
Hence, quotient the $\mathrm{SU}(2)$ action on $\mathcal{N}$, the fiber is isomorphic to
\[H^1\big(\pi_1(\Sigma_g),\mathfrak{h}_{\mathrm{Ad}_{\rho_0}}^{\perp}\big)\ \cong\ \frac{Z^1\big(\pi_1(\Sigma_g),\mathfrak{h}_{\mathrm{Ad}_{\rho_0}}^{\perp}\big)}{B^1\big(\pi_1(\Sigma_g),\mathfrak{h}_{\mathrm{Ad}_{\rho_0}}^{\perp}\big)}.\] 
On the fiber $\mathcal{N}_{\rho_0}$, the pre-symplectic form $\omega_{\rho_0}^N$ determines a short exact sequence of spaces,
\[B^1\big(\pi_1(\Sigma_g),\mathfrak{h}_{\mathrm{Ad}_{\rho_0}}^{\perp}\big)\ \rightarrow\ Z^1\big(\pi_1(\Sigma_g),\mathfrak{h}_{\mathrm{Ad}_{\rho_0}}^{\perp}\big)\ \xrightarrow{\omega_{\rho_0}^N} H^1\big(\pi_1(\Sigma_g),\mathfrak{h}_{\mathrm{Ad}_{\rho_0}}^{\perp}\big),\]
where 
\[\ker{\omega_{\rho_0}^N}\ =\ B^1\big(\pi_1(\Sigma_g),\mathfrak{h}_{\mathrm{Ad}_{\rho_0}}^{\perp}\big).\]
Hence, the vector space $H^1\big(\pi_1(\Sigma_g),\mathfrak{h}_{\mathrm{Ad}_{\rho_0}}^{\perp}\big)$ is symplectic with the symplectic form $\Omega_{\rho_0}^N$. Therefore, the bundle 
\[\mathcal{N}/\mathrm{SU}(2)\ \rightarrow\ \mathcal{K}_0\]
is a symplectic vector bundle with the canonical Narasimhan-Atiyah-Bott-Goldman symplectic form.

\end{proof}

\begin{theorem}
\label{SymplecticNormalConeInCharacterVariety}
    The symplectic vector space $\Big(H^1\big(\pi_1(\Sigma_g),\mathfrak{h}_{\mathrm{Ad}_{\rho_0}}^{\perp}\big),\Omega_{[\rho_0]}^N\Big)$ is a Hamiltonian $\mathrm{U}(1)$-space, where the moment map is
    \[\begin{aligned}
        \Phi_{[\rho_0]}:\ H^1\big(\pi_1(\Sigma_g),\mathfrak{h}_{\mathrm{Ad}_{\rho_0}}^{\perp}\big)\ &\rightarrow\ \mathfrak{h}^*\\
        \tau\ &\mapsto\ [\tau,\tau].
    \end{aligned}\]
\end{theorem}
\begin{proof}
The symplectic vector space is 
\[\Big(H^1\big(\pi_1(\Sigma_g),\mathfrak{h}_{\mathrm{Ad}_{\rho_0}}^{\perp}\big),\Omega_{[\rho_0]}^N\Big).\]
Let $\xi\in\mathfrak{h}$, and $\tau\in H^1\big(\pi_1(\Sigma_g),\mathfrak{h}_{\mathrm{Ad}_{\rho_0}}^{\perp}\big)$, then the moment map $\Phi_{[\rho_0]}$ satisfies the equation
\[\langle\Phi_{[\rho_0]}(v),\xi\rangle\ =\ \frac{1}{2}\Omega_{[\rho_0]}^N(\mathrm{ad}_{\xi}\tau,\tau).\]
In particular, the operator $\mathrm{ad}_{\xi}$ is the 0-degree adjoint transformation, and
\[\mathrm{ad}_{\xi}(u\wedge u)\ =\ (\mathrm{ad}_{\xi}u)\wedge u-u\wedge(\mathrm{ad}_{\xi}u)\ =\ 2(\mathrm{ad}_{\xi}u)\wedge u.\]
Hence,

\[\begin{aligned}
    \Omega_{[\rho_0]}^N(\mathrm{ad}_{\xi}\tau,\tau)\ &=\ \Omega_{\rho_0}^N([\xi,\tau],\tau)\ =\ \frac{1}{2}\Omega_{\rho_0}^N([\tau,\tau], \xi).
\end{aligned}\]
Indeed, 
\[H^2\big(\pi_1(\Sigma_g),\mathfrak{h}\big)\ \cong\ \mathfrak{h}^*.\]
Therefore, using the Fox derivative, up to a scalar, the moment map is
\[\Phi_{[\rho_]0}(\tau)\ :=\ [\tau,\tau].\]
In parcular, $\Phi_{[\rho_0]}$ is $\mathrm{U}(1)$-equivariant (\ref{EquivariantIsotropyGroupActionOnCone}).
\end{proof}
Apply the Marsden-Weinstein-Meyer theorem to the Hamiltonian space 
   \[\Big(H^1\big(\pi_1(\Sigma_g),\mathfrak{h}_{\mathrm{Ad}_{\rho_0}}^{\perp}\big), \Omega_{[\rho_0]}^N, \mathrm{Ad}_{\mathrm{U}(1)},\Phi_{[\rho_0]}\Big),\]
   the symplectic reduction determines the symplectic structure of the reduced space.
\begin{theorem}
\label{SymplecticIsotropyGroupReduction}
Let $i:\Phi_{[\rho_0]}^{-1}(0)\hookrightarrow\ H^1\big(\pi_1(\Sigma_g),\mathfrak{h}_{\mathrm{Ad}_{\rho_0}}^{\perp}\big)$ be the inclusion map, where $\mathrm{U}(1)$ acts freely. Then the reduced space is symplectic:
\[\big(\Phi_{[\rho_0]}^{-1}(0)/\mathrm{U}(1),\widetilde{\Omega}^N_{[\rho_0]}\big),\]
and $i^*\Omega_{\rho_0}^N=\pi ^*\widetilde{\Omega}_{[\rho_0]}^N$, where $\pi:\Phi_{[\rho_0]}^{-1}(0)\rightarrow \Phi_{[\rho_0]}^{-1}(0)/\mathrm{U}(1)$. In particular,
    \[\dim_{\mathbb{R}}\Phi_{[\rho_0]}^{-1}(0)/\mathrm{U}(1) =\ 4g-6.\]
\end{theorem}
\begin{proof}
    The space of 1-coboundaries $B^1\big(\pi_1(\Sigma_g),\mathfrak{h}_{\mathrm{Ad}_{\rho_0}}^{\perp}\big)$ is a subspace of the quadratic cone $\mathcal{C}_{\rho_0}$. The set of classes of 1-cocycles in $H^1\big(\pi_1(\Sigma_g),\mathfrak{h}_{\mathrm{Ad}_{\rho_0}}^{\perp}\big)$ is denoted by $\mathcal{C}'_{[\rho_0]}$, and $\mathcal{C}'_{[\rho_0]}\ =\ \Phi_{[\rho_0]}^{-1}(0)$. The symplectic reduction is
   \[q:\ \mathcal{C}'_{[\rho_0]}\ \rightarrow\ \widetilde{\mathcal{C}}_{[\rho_0]},\]
   where $\widetilde{\mathcal{C}}_{[\rho_0]}\ :=\ \mathcal{C}'_{[\rho_0]}/\mathrm{U}(1)$.   Hence, $\widetilde{\mathcal{C}}_{[\rho_0]}$ is symplectic with the symplectic form $\widetilde{\Omega}^N_{[\rho_0]}$ so that
\[q^{*}\widetilde{\Omega}^N_{[\rho_0]}\ =\ i^*\Omega_{[\rho_0]}^N.\]
   \end{proof}

\begin{theorem}
\label{A Reduced Symplectic Cone Bundle Over The Kummer Variety}
    \textit{The reduced cone bundle 
    \[\widetilde{C}\ \rightarrow\ \mathcal{K}_0\ \]
    is a symplectic with the symplectic 2-form $\widetilde{\Omega}^{N}$, where the fiber at $[\rho_0]$ is the symplectic quadratic cone
    \[\big(\widetilde{\mathcal{C}}_{[\rho_0]}, \widetilde{\Omega}^N_{[\rho_0]}\big).\]}
\end{theorem}

\begin{proof}
 Let $\rho\in\mathcal{B}_{\varepsilon}(\rho_0)$ that an open ball of $\rho_0$ in $\mathcal{Y}$ so that
\[\rho\ =\ \exp{\big(\varepsilon u+\mathcal{O}(\varepsilon^2)\big)}\cdot\rho_0.\]
Suppose $u_1$ is in $Z^1\big(\pi_1(\Sigma_g),\mathfrak{h}^{\perp}_{\mathrm{Ad}_{\rho_0}}\big)$ so that $[u_1,u_1]=0$ in $H^2\big(\pi_1(\Sigma_g),\mathfrak{h}\big)$. Under the local trivialization map $R_{\exp{(\varepsilon u)}}$,
\[[R_{\exp{(\varepsilon u)}}u_1, R_{\exp{(\varepsilon u)}}u_1]\ =\ 0.\]
Hence, the map $R_{\exp{(\varepsilon u)}}$ preserves the quadratic bracket. On the other hand, the quadratic bracket preserves the transition function:
\[[\ \cdot\wedge\cdot\ ]:\ \pi_{C}^{-1}\big(B_{\varepsilon}(\rho_0)\big)\ \rightarrow\ \pi_{C}^{-1}\big(B_{\varepsilon}(\rho_0)\big).\]
Furthermore, the reduction is compatible with the bundle structure, 
\[R_{\exp{\varepsilon u}}[\mathrm{Ad}_{h}u_1,\mathrm{Ad}_{h}u_1]\ =\ \mathrm{Ad}_h[R_{\exp{\varepsilon u}}u_1,R_{\exp{\varepsilon u}}u_1].\]
Therefore, the bundle $\mathcal{C}'\rightarrow\mathcal{K}_0$ is a subbundle of the quotient bundle $\mathcal{N}/\mathrm{SU}(2)$,
\[\begin{tikzcd}
	{\widetilde{C}} & {\mathcal{N}/\mathrm{SU}(2)} \\
	{\mathcal{K}} & {\mathcal{K}}
	\arrow[from=1-1, to=2-1]
	\arrow[from=1-2, to=2-2]
	\arrow["{\mathrm{Id}}", from=2-1, to=2-2]
\end{tikzcd}\]
which covers the identity on the base.

\end{proof}

\subsection{The symplectic form on the stratum $\chi_0$}
Let $\rho_0$ be in $\mathcal{Y}_0$. The locus $\mathcal Y_0$ carries a Weyl group symmetry. The projection
\[\mathcal Y_0\longrightarrow \chi_0\]
is governed by this Weyl action and is a $W$-principal fiberation (\ref{WeylGroupDefinition}). The map between the reductions over $\chi_0$ is a map of bundles
\[\begin{tikzcd}
	{\mathcal{Y}_0} & {\widetilde{\mathcal{Y}}_0} \\
	& {\chi_0}
	\arrow[hook, from=1-1, to=1-2]
	\arrow[from=1-1, to=2-2]
	\arrow[from=1-2, to=2-2]
\end{tikzcd}\]
The right bundle has the fiber $\mathrm{SU}(2,A_2)^0/\mathrm{U(1)}$, where $A_2=\mathbb{R}[[t]]/(t^3)$, $\mathfrak{m}_2=(t,t^2)$ and
\[\mathrm{SU}(2,A_2)^0\ =\ \exp{\big(\mathfrak{su}(2)\otimes \mathfrak{m}_2\big)}.\]
\begin{proposition}
\label{Compatibility}
    \textit{The q-Hamiltonian space $\big(\mathcal{Y}_0, \mathrm{Ad}_{W}, \omega, \mu|_{\mathcal{Y}_0}\big)$ sits in the tubular neighborhood $\mathcal{U}_{\mathcal{R}}$.
     In particular, $\omega^{\mathcal{Y}_0}$ is reduced to $\mathcal{Y}_0/W$ under the q-Hamiltonian reduction, $q^*\widetilde{\Omega}_{[\rho_0]}=\omega^{\mathcal{Y}_0}$ where
     \[q:\ \mathcal{Y}_0\ \rightarrow\ \chi_0.\] }
\end{proposition}
\begin{proof}
\begin{itemize}
\item Indeed,
\[\mathcal{Y}_0\ =\  \mathrm{Hom}\big(\pi_1(\Sigma_g),\mathrm{SU}(2)\big)_{\mathrm{U}(1)},\]
and $\mu(\mathcal{Y}_0)\equiv1$. Moreover,
\[d\mu\big|_{\mathcal{Y}_0}\ \equiv\ 0.\]
Meanwhile, the fundamental vector fields of $\mathbb{Z}_2$-action are all $0$'s. 
\item Secondly, $\omega$ is closed,
\[d\omega\ =\ 0.\]
Meanwhile, $\mu^*\big|_{\mathcal{Y}_0}=0$. Hence,
\[d\omega|_{\mathcal{Y}_0}\ =\ -\mu^{*}\chi,\]
where $\chi$ is the canonical 3-form 
\[\chi\ =\ \frac{1}{12}(\theta,[\theta,\theta])\ =\ \frac{1}{12}(\bar{\theta},[\bar{\theta},\bar{\theta}]).\]

\item Along $\mathcal{Y}_0$, the tangent space is
\[T_{\rho_0}\mathcal{Y}\ \cong\ Z^1(\Gamma_g, \mathfrak{h}).\]
The solution the equation
\[\mathrm{Ad}_{\mu(\mathcal{Y})}\xi+\xi\ =\ 0. \]
is $\xi=0$. At the same time, $B^1(\Gamma_g,\mathfrak{h})=\{0\}$.
Therefore,
\[\ \ker{\omega}|_{\mathcal{Y}_0}\ =\ \{0\}\ =\ B^1(\Gamma_g,\mathfrak{h}) .\]
\end{itemize}
Therefore, $\big(\mathcal{Y}_0, \mathrm{Ad}_{w}, \omega, \mu\big)$ is a q-Hamiltonian space with the standard q-Hamiltonian reduction.
\end{proof}

\begin{theorem}
\label{SymplecticEmbedding}
    The symplectic cone bundle $\widetilde{\mathcal{C}}\rightarrow\mathcal{K}_0$ sits in the character variety $\chi\big(\pi_1(\Sigma_g),\mathrm{SU}(2)\big)$, where $\mathcal{K}_0=\chi_0$. Therefore, the symplectic structure of $\chi_0$ is
    \[\widetilde{\Omega}_{\chi_0}\ =\ \omega^{\mathcal{K}_0}+\widetilde{\Omega}^N.\]
    \end{theorem}
    
\begin{remark}
    The following diagram visualizes the whole picture of this paper.
    \[\begin{tikzcd}
	&& {\widetilde{\mathcal{Y}}_0} & {\mathsf{Hom}\big(\pi_1(\Sigma_g),\mathrm{SU}(2)\big)} \\
	{\mathcal{N}} & {\mathcal{Y}_0} & {\chi_0} \\
	{\mathcal{N}/\mathrm{SU}(2)} & {\mathcal{K}_0}\\
    {\widetilde{\mathcal{C}}} &{\mathcal{K}_0}
	\arrow[hook, from=1-3, to=1-4]
	\arrow["{/\mathrm{SU}(2)}", from=1-3, to=2-3]
	\arrow[from=2-1, to=2-2]
	\arrow[hook, from=2-2, to=1-3]
	\arrow["/W", from=2-2, to=2-3]
	\arrow["{/\mathrm{SU}(2)}"', from=2-2, to=3-2]
	\arrow[ from=3-1, to=3-2]
	\arrow[from=3-2, to=2-3]
    \arrow["{\mathrm{Id}}",from=4-2, to=3-2]
    \arrow[from=4-1,to=4-2]
\end{tikzcd}\]
The embedding of bundles enables the symplectic form to be the sum of tangential and normal directions.
\end{remark}

\newpage
\printbibliography

@book{FranciosRepresentation,
author = {Labourie, François},
year = {2013},
month = {12},
pages = {},
title = {Lectures on Representations of Surface Groups},
isbn = {978-3-03719-127-9},
doi = {10.4171/127}
}

@article{BillMillsonDeforamtion,
author = {Goldman William M., Millson John J.},
journal = {Publications Math\'ematiques de l'IH\'ES},
language = {eng},
pages = {43-96},
publisher = {Institut des Hautes \'Etudes Scientifiques},
title = {The deformation theory of representations of fundamental groups of compact Kähler manifolds},
url = {http://eudml.org/doc/104031},
volume = {67},
year = {1988},
}

@article{Alekseev1997LieGV,
  title={Lie group valued moment maps},
  author={Anton Yu . Alekseev and Anton Malkin and Eckhard Meinrenken},
  journal={Journal of Differential Geometry},
  year={1997},
  volume={48},
  pages={445-495},
  url={https://api.semanticscholar.org/CorpusID:18008533}
}

@book{KobayashiDifferentialGeometry1,
  title={Foundations of Differential Geometry, Vol. I},
  author={Kobayashi, Shoshichi and Nomizu, Katsumi},
  year={1963},
  publisher={Interscience Publishers, New York-London},
  series={Interscience Tracts in Pure and Applied Mathematics},
  number={15}
}

@article{SjamaarStratified,
  title={Stratified symplectic spaces and reduction},
  author={Sjamaar, Reyer and Lerman, Eugene},
  journal={Annals of Mathematics},
  pages={375--422},
  year={1991},
  publisher={JSTOR}
}

@article{SikoraCharacterVarieties,
  title={Character Varieties},
  author={Adam S. Sikora},
  journal={arXiv preprint arXiv:0902.2589},
  year={2009},
  url={https://arxiv.org/abs/0902.2589}
}

@article{BillSymplecticNature,
  title={The Symplectic Nature of Fundamental Groups of Surfaces},
  author={William M. Goldman},
  journal={Advances in Mathematics},
  year={1984},
  volume={54},
  pages={200-225},
  doi={10.1016/0001-8708(84)90040-9},
  url={https://api.semanticscholar.org/CorpusID:16484957}
}

@article{JeffreyExtendedModuliSpace,
  author    = {Lisa C. Jeffrey},
  title     = {Symplectic Forms on Moduli Spaces of Flat Connections on 2-manifolds},
  journal   = {arXiv preprint alg-geom/9404013},
  year      = {1994},
  url       = {https://arxiv.org/abs/alg-geom/9404013}
}

@InProceedings{HuebschmannSingularetiesPossion,
author="Huebschmann, Johannes",
editor="Landsman, N. P.
and Pflaum, M.
and Schlichenmaier, M.",
title="Singularities and Poisson geometry of certain representation spaces",
booktitle="Quantization of Singular Symplectic Quotients",
year="2001",
publisher="Birkh{\"a}user Basel",
address="Basel",
pages="119--135",
isbn="978-3-0348-8364-1"
}

@article{HuebschmannPoissonSU2,
  author  = {Huebschmann, Johannes},
  title   = {Poisson Geometry of Flat Connections for {$\mathrm{SU}(2)$}-Bundles on Surfaces},
  journal = {Mathematische Zeitschrift},
  volume  = {221},
  number  = {2},
  pages   = {243--260},
  year    = {1996}
}

@article{Artin1968,
  author  = {M. Artin},
  title   = {On the Solutions of Analytic Equations},
  journal = {Inventiones mathematicae},
  volume  = {5},
  number  = {4},
  pages   = {277--291},
  year    = {1968},
  url     = {https://eudml.org/doc/141922}
}
\end{document}